\documentclass{amsart}  \usepackage{amssymb,latexsym}
\usepackage{graphicx}
\usepackage{etex, verbatim}
\usepackage{color}
\usepackage[pagewise]{lineno}
\usepackage{perpage} 
\usepackage[margin=1in]{geometry}
\usepackage{array}
\usepackage[all,2cell,ps]{xy}
\newcommand{\edge}[1]{\ar@{-}[#1]}

\newcommand{\lr}{\mbox{$\longrightarrow$}}

\newcommand{\cf}{{\mathcal F}}

\newcommand{\lra}{\longrightarrow}
\newcommand{\da}{\Big\downarrow}
\newcommand{\ra}{ \rightarrow}

\newcommand{\hra}{\hookrightarrow}

\newcommand{\be}{\begin{equation}}
\newcommand{\ee}{\end{equation}}
\newcommand{\K}{{\mbox{K}}}

\newcommand{\R}{\mbox{R}}
\newtheorem{guess}{Theorem}[section]
\newcommand{\bth}{\begin{guess}$\!\!\!${\bf }~}
\newcommand{\eeth}{\end{guess}}
\renewcommand{\bar}{\overline}
\newtheorem{propo}[guess]{Proposition}
\newcommand{\bpropo}{\begin{propo}$\!\!\!${\bf }~}
\newcommand{\epropo}{\end{propo}}

\newtheorem{lema}[guess]{Lemma}
\newcommand{\blem}{\begin{lema}$\!\!\!${\bf }~}
\newcommand{\elem}{\end{lema}}

\newtheorem{defe}[guess]{Definition}
\newcommand{\bdefe}{\begin{defe}$\!\!\!${\bf }~}
\newcommand{\edefe}{\end{defe}}

\newtheorem{notation}[guess]{\bf Notation}

\newtheorem{coro}[guess]{Corollary}
\newcommand{\bcor}{\begin{coro}$\!\!\!${\bf }~}
\newcommand{\ecor}{\end{coro}}

\newtheorem{rema}[guess]{Remark}
\newcommand{\brem}{\begin{rema}$\!\!\!${\bf }~\rm}
\newcommand{\erem}{\end{rema}}

\newtheorem{notn}[guess]{Notation}
\newcommand{\bnotn}{\begin{notn}$\!\!\!${\bf }~\rm}
\newcommand{\enotn}{\end{notn}}

\newtheorem{exam}[guess]{Example}
\newcommand{\beg}{\begin{exam}$\!\!\!${\bf }~\rm}
\newcommand{\eeg}{\end{exam}}

\newcommand{\cb}{{\mathcal B}}

\newcommand{\bz}{\mathbb{Z}}

\newcommand{\cv}{{\mathcal V}}
\newcommand{\ce}{{\mathcal E}}
\newcommand{\ck}{{\mathcal K}}

\newcommand{\cx}{{\mathcal X}}
\newcommand{\cz}{{\mathcal Z}}

\newcommand{\cg}{{\mathcal G}}
\newcommand{\cw}{\mathcal{W}}
\newcommand{\ci}{{\mathcal I}}
\newcommand{\cy}{{\mathcal Y}}

\newcommand{\br}{\mathbb{R}}
\newcommand{\tG}{\widetilde{G}}
\newcommand{\tB}{\widetilde{B}}
\newcommand{\tC}{\widetilde{C}}
\newcommand{\tT}{\widetilde{T}}

\newcommand{\tW}{\widetilde{W}}

\newcommand{\tphi}{\widetilde{\Phi}}

\newcommand{\bc}{\mathbb{C}}
\newcommand{\bp}{\mathbb{P}}

\renewcommand{\phi}{\varphi}
\newcommand{\co}{{\mathcal O}}

\begin{document}
\title{$\K$-theory of regular compactification bundles}
\author{V. Uma}
\address{Department of Mathematics, Indian Institute of Technology-Madras, Chennai, India}
\email{vuma@iitm.ac.in}

\begin{abstract}
  Let $G$ be a split connected reductive algebraic group. Let
  $\ce\lra \cb$ be a $G\times G$-torsor over a smooth base scheme
  $\cb$ and $X$ be a regular compactification of $G$. We describe the
  Grothendieck ring of the associated fibre bundle
  $\ce(X):=\ce\times_{G\times G} X$, as an algebra over the
  Grothendieck ring of a canonical toric bundle over a flag bundle
  over $\cb$. These are relative versions of the results in
  \cite{u1,u2}, and generalize the classical results on the
  Grothendieck rings of projective bundles, toric bundles \cite{su}
  and flag bundles \cite{fl,pr}.
\end{abstract}

\maketitle 

\thispagestyle{empty}

\footnote{{\bf
keywords:} Regular compactification bundles, toric bundles,
flag bundles, $K$-theory}

\section*{Introduction}
In this article we shall consider algebraic groups and varieties to be
defined over an algebraically closed field $F$. By an algebraic
group we mean an {\it affine, smooth group scheme of finite
  type} over $F$.  By a scheme we mean a separated scheme of finite
type over $F$. By a variety we mean an integral scheme of finite type
over $F$. All schemes we consider are assumed to be noetherian. By a
{\it torsor} we shall always mean {\it locally isotrivial} in the
sense of \cite{serre}.

Let $G$ denote a connected reductive algebraic group. Let $C$ be the
center of $G$ and let $G_{ad}:=G/C$ be the corresponding semisimple
adjoint group.

A normal complete variety $X$ is called an {\it equivariant
compactification} of $G$ if $X$ contains $G$ as an open subvariety and
the action of $G\times G$ on $G$ by left and right multiplication
extends to $X$.  We say that $X$ is a {\it regular compactification}
of $G$ if $X$ is an equivariant compactification of $G$ which is
regular as a $G\times G$-variety ( \cite[Section 2.1]{Br}). Smooth
complete toric varieties are regular compactifications of the
torus. For the adjoint group $G_{ad}$, the wonderful compactification
$\bar {G_{ad}}$ constructed by De Concini and Procesi in \cite{DP1} is
the unique regular compactification of $G_{ad}$ with a unique closed
$G_{ad}\times G_{ad}$-orbit.

Let $\mathcal{E}\lra \mathcal{B}$ be a $G\times G$-torsor over a base
scheme $\cb$. Let $X$ be a projective regular compactification of a
connected reductive algebraic group $G$. Let
$\ce(X):=\ce\times_{(G\times G)} X$ denote the associated bundle with
fibre $X$ and base $\cb$. Since $\ce$ is the total space of a
$G\times G$-torsor over $\cb$, it is a $G\times G$-scheme. Further,
the space $\ce(X)$ also gets the structure of a scheme (see
\cite[Proposition 23]{eg}). Moreover, since $X$ is a smooth
$G$-scheme, when the base scheme $\cb$ is smooth, it can be seen that
$\ce(X)$ is a smooth scheme (see below for details).

The main aim of this article is to describe the Grothendieck ring of
algebraic vector bundles on $\ce(X)$ as an algebra over the
Grothendieck ring of algebraic vector bundles on a smooth base scheme
$\cb$.  This is with a view to generalize and is motivated by the
corresponding classical results on projective bundles \cite{BGI},
toric bundles \cite{su}, and flag bundles \cite{fl,pr}. This is also a
relative version of the results in \cite{u1,u2} on $\K$-theory of
regular compactifications when $\cb=pt$.

We now fix some notations before we give an overview of the
main results.

For a $G$-scheme $X$, let $\K^0_{G}(X)$ denote the Grothendieck ring
of $G$-equivariant vector bundles on $X$ and $\K_0^G(X)$ denote the
Grothendieck group of $G$-equivariant coherent sheaves on $X$. Let
$\R(G)$ denote the Grothendieck ring of representations of $G$ over
the field $F$. Identifying $\R(G)$ with $\K_{G}^0(pt)$, the pullback
by the projection $X\lra pt$ gives a canonical map
$\R(G)\lra \K^0_{G}(X)$ . This gives $\K^0_{G}(X)$ the structure of an
$\R(G)$-algebra. Also $\K_0^{G}(X)$ is a module over the ring
$\K^0_{G}(X)$ and is hence an $\R(G)$-module. Moreover, when $X$
is smooth $\K^0_{G}(X)\cong \K_0^G(X)$.

In Section 1 we prove some preliminary results on the $G$-equivariant
$\K$-theory of schemes having equivariant cellular structure, assuming
that $\pi_1(G)$ is torsion free. Let $X$ be any $G$-cellular scheme,
that is a $G$-scheme with paving by $G$-stable affine cells (see
Definition \ref{gcell}).

We prove an equivariant Kunneth theorem (see Theorem
\ref{equivkunneth}) where we show that we have an isomorphism
\(\displaystyle \K^G_0(Y)\bigotimes_{R(G)} \K_0^{G}(X)
\cong\K_0^G(Y\times X) \) of $\K_G^0(Y)$-modules. Moreover, when $X$
and $Y$ are smooth then the Kunneth map is an isomorphism of
$\K_{G}^0(Y)$-algebras.

Let $p:\ce\lra\cb$ denote a $G$-torsor over a base scheme $\cb$ and
$X$ be a $G$-scheme. Then by \cite[Proposition 23]{eg}, $\ce$ is a
$G$-scheme and the associated bundle $\ce(X):=\ce\times_{G} X$ is a
scheme. We have the projection $\pi:\ce(X)\lra \cb$ defined as
$\pi[e,x]=p(e)$ for $[e,x]\in \ce(X)$. Furthermore, when $\cb$ and $X$
are smooth, $\ce(X)$ is a smooth scheme. (Since $G$ is smooth, the
quotient map $\ce\lra \cb$ is a smooth morphism. Thus when the base
scheme $\cb$ is smooth, by the local isotriviality of the $G$-torsor
$\ce\lra \cb$, it follows that the total space $\ce$ is smooth. When
$X$ is a smooth $G$-scheme, $\ce\times X$ is smooth. Thus the quotient
$\ce(X)$ of $\ce\times X$ by the free action of $G$ is a smooth
scheme.)

For $X$ a $G$-cellular scheme (resp. a $G$-scheme with a $T$-cellular
structure, where $T$ denotes a maximal torus of $G$ acting on $X$ by
restriction), in Corollary \ref{maincor} (resp. Corollary
\ref{corollary2}) we derive the structure of the Grothendieck group
$\K_0(\ce(X))$ of coherent sheaves on $\ce(X)$, as a module over the
Grothendieck ring $\K^0(\cb)$ of algebraic vector bundles on $\cb$. We
further note that Corollary \ref{corollary2} holds for any projective
$G$-scheme $X$ having finitely many $T$-fixed points.

In particular, when $X$ is smooth and the base scheme $\cb$ is smooth
we get the structure of $\K_0(\ce(X))=\K^0(\ce(X))$ as a
$\K^0(\cb)$-algebra.  More precisely, using the Kunneth formula and
the fact that $\K^0_{G}(\ce)=\K^0(\ce/G)=\K^0(\cb)$, we show that
$\K^0(\ce(X))$ is isomorphic to a canonical extension of scalars of
the $\R(G)$-algebra $\K^0_{G}(X)$, to the ring $\K^0(\cb)$.

Let $\cb$ be a smooth scheme and let $p:\ce\lra \cb$ be a $T$-torsor
and $X$ be a smooth projective $T$-variety with finitely many
$T$-fixed points and invariant curves. Let $X^T$ denote the set of
$T$-fixed points in $X$. If $|X^T|=m$, we consider
\(\displaystyle\prod_{i=1}^m \K^0(\cb)\) which is a ring under pointwise
addition and multiplication. It further gets a canonical
$\K^0(\cb)$-algebra structure through the diagonal inclusion. The ring
$\K^0(\ce(X))$ also has a $\K^0(\cb)$-algebra structure induced by pull
back of vector bundles under the projection $\pi$. We next prove
Theorem \ref{relloc} which is a relative version of the localization
theorem for $\K^0_{T}(X)$ (see \cite{VV}, \cite[Theorem 1.3]{u1}). Here
we show that the restriction to the $T$-fixed points of the fiber
induces a canonical inclusion
\[ \K^0(\ce(X)\hra \prod_{i=1}^m \K^0(\cb)\] of $\K(\cb)$-algebras. We
further show that the image of $\K^0(\ce(X))$ in
\(\displaystyle\prod_{i=1}^m \K^0(\cb)\) is precisely the intersection of
the images of the Grothendieck rings of $\bp^1$-bundles associated to
$\ce\stackrel{p}{\lra} \cb$, corresponding to the projective lines
joining two distinct $T$-fixed points.

In Section 2 we prove our main results. Let $G$ be a connected
reductive algebraic group and let $\tG$ be a factorial cover of
$G$ (see \ref{fc}). In particular, $\pi_1(\tG)$ is torsion free.

Let $p:\ce\lra \cb$ a $\tG\times \tG$-torsor
over a smooth base scheme $\cb$. We let $X$ to be a projective regular
compactification of $G$. Let $\ce(X):=\ce\times_{\tG\times \tG} X$
where $\tG$ acts on $X$ via its canonical projection to $G$.

Let $T$ denote a maximal torus of $G$ and $B$ a Borel subgroup
containing $T$. Let $W$ denote the Weyl group of $(G,T)$.  In Theorem
\ref{first}, using Theorem \ref{equivkunneth} and \cite[Corollary
2.3]{u1} we describe the Grothendieck ring of $\ce(X)$ as
$\mbox{diag}(W)$-invariants of the Grothendieck ring of a toric
bundle, with fibre the toric variety $\bar{T}\subseteq \bar{G}=X$, and
base another bundle over $\cb$ with fibre $G/B^-\times G/B$. We note
that here the $\mbox{diag}(W)$-action on the Grothendieck ring of the
toric bundle is induced from its canonical action on $\bar{T}$ (see
\cite[Proposition A1, A2 ]{Br}). This is the relative version of
\cite[Proposition 2.15]{u1}.

In Theorem \ref{second}, we use Theorem \ref{equivkunneth}, Theorem
\ref{relloc}, \cite[Corollary 2.  2]{u1} and \cite[Theorem 2.4]{u2},
to further describe the multiplicative structure of $\K^0(\ce(X))$, as
an algebra over the Grothendieck ring of a toric bundle with fibre the
toric variety $\bar{T}^+$, and base a flag bundle. The toric variety
$\bar{T}^+$ is associated to a smooth fan in the lattice of one
parameter subgroups of $T$, supported on the positive Weyl chamber
(see \cite[Proposition A1, A2]{Br}).  This is the relative version of
\cite[Theorem 3.1]{u2}. Finally, in Section 3 we retrieve the known
results on the Grothendieck ring of toric bundles and flag bundles.

\noindent {\bf Acknowledgements:} The author is grateful to
Prof. Michel Brion for his patient reading and invaluable comments and
suggestions for improvement of the earlier versions of this
manuscript. The author is grateful to Prof. D. S. Nagaraj for reading
the manuscript and for some valuable comments and suggestions. The
author is grateful to Prof. V. Balaji for several valuable
discussions. The author is especially grateful to the referee for a very
careful reading of the manuscript and for several invaluable comments
and suggestions for improvement.

\section{Preliminaries on equivariant $\K$-theory}

Throughout this section we shall assume that $G$ is a connected
reductive algebraic group. In Subsections 1.1 and 1.2 we shall assume in
addition that $\pi_1(G)$ is torsion free. This is with a view to apply
the results in \cite{mer} which require this hypothesis.

\subsection{Kunneth formula for equivariant $\K$-theory}

For $G$-schemes $X$ and $Y$, the map
\be\label{etp1}\boxtimes: \K_0^{G}(Y)\bigotimes_{\bz}\K_0^{G}(X)\lra \K_0^{G}(Y\times X).\ee
induced by the external tensor product of $G$-equivariant coherent sheaves
is defined by
\be\label{etp2}(\cg,\cg')\mapsto \cg\boxtimes
  \cg':=p_Y^*(\cg)\bigotimes_{\co_{Y\times X}}p_X^*(\cg')\ee where
$p_X$ and $p_Y$ are the projections from $Y\times X$ to $X$ and $Y$
respectively. Moreover, since $p^*_X$ and $p^*_Y$ are $\R(G)$-module
maps,  the elements of the form $1\bigotimes a-a\bigotimes 1$
for $a\in \R(G)$ map to $0$ under $\boxtimes$. This induces a map
of $\R(G)$-modules\be\label{extalg}\varphi:
\K_0^{G}(Y)\bigotimes_{\R(G)} \K_0^{G}(X) \lra \K_0^{G}(Y\times
X).\ee

\bdefe\label{gcell}  A {\em
  $G$-cellular scheme} is a $G$-scheme $X$ equipped with
a $G$-stable algebraic cell decomposition. In other words
there is a filtration
\[X=X_1\supseteq X_2\supseteq\cdots\supseteq X_m\supseteq X_{m+1}=\emptyset\] where
each $X_i$ is a closed $G$-stable subscheme of $X$ and
$X_i\setminus X_{i+1}=U_i$ is $G$-equivariantly isomorphic to the
affine space $\mathbb{A}_{F}^{k_i}$ equipped with a linear action of
$G$, for $1\leq i\leq m$. \edefe

For $Y$ any $G$-scheme and for $X_i$ and $U_i$
for $1\leq i\leq m$, as in Definition \ref{gcell}, we have the
following maps of $\R(G)$-modules
\be\label{map}\phi_i:\K^{G}_0(Y)\bigotimes_{\R(G)}
\K^{G}_0(X_i)\lra \K_0^{G}(Y\times X_i)\ee and \be
\psi_i:\K_0^{G}(Y)\bigotimes_{\R(G)} \K_0^{G}(U_i)\lra
\K_0^{G}(Y\times U_i).\ee

We now show that the $\R(G)$-modules on either side of (\ref{extalg})
have a canonical structure of $\K_{G}^0(Y)$-modules. This is induced
from the $\K_G^0(Y)$-module structure on $\K^G_0(Y)$.

For every $l\geq 0$, we have
$p_{Y}^*: \K^{G}_l(Y)\lra \K^{G}_l(Y\times X)$ induced by the pull
back of $G$-equivariant coherent sheaves under $p_{Y}$. In particular,
$p_{Y}^*: \K^{G}_0(Y)\lra \K^{G}_0(Y\times X)$ maps $\K^0_{G}(Y)$ to
$\K^0_{G}(Y\times X)$ via pull back of $G$-equivariant vector bundles
under $p_{Y}$ giving $\K_{G}^0(Y\times X)$ the structure of a
$\K^0_G(Y)$-algebra.

Furthermore, since $\K_l^{G}(Y\times X)$ is a
$\K_{G}^0(Y\times X)$-module (induced by tensor product of
$G$-equivariant vector bundles with $G$-equivariant coherent sheaves
see \cite[p.247-248]{CG}, \cite{Q}), it follows that
$\K_l^{G}(Y\times X)$ gets the structure of a
$\K_{G}^0(Y)$-module for $l\geq 0$. Furthermore, $p_{Y}^*$ is a
morphism of $\K_{G}^0(Y)$-modules.

On the other hand we note that
$\K_0^{G}(Y)\bigotimes_{\R(G)} \K_0^{G}(X)$ has a canonical
structure of
$\K^0_{G}(Y)\bigotimes_{\R(G)} \R(G)\cong \K^0_{G}(Y)$-module.

We now claim that $\phi$ preserves the $\K_{G}^0(Y)$-module
structure. Note that for a $G$-equivariant vector bundle $E$ on
$Y$, $p_{Y}^*(E)$ is a $G$-equivariant vector bundle on $Y\times X$
and
$p_{Y}^*(E\otimes_{\co_{Y}} \cf)=p_{Y}^*(E)\otimes_{\co_{X\times
    Y}}p_{Y}^*(\cf)$ for every $G$-equivariant coherent sheaf $\cf$ on
$Y$. The claim follows from the $\K_{G}^0(Y)$-module structure on
$\K^G_0(Y)\bigotimes_{\R(G)} \K^G_0(X)$ and $\K^{G}_0(Y\times X)$.

We recall below the Thom isomorphism theorem in higher $G$-equivariant
$K$-theory (see \cite{Q}, \cite{Th} or \cite[Theorem 5.4.17]{CG}).

{\bf Theorem:} {\it (Thom isomorphism)} Let $\pi':E\longrightarrow X$
be a $G$-equivariant affine bundle which is a torsor under the vector
group scheme associated with a $G$-equivariant vector bundle over a
$G$-scheme $X$. For any $j\geq 0$ the morphism
${\pi'}^*:\K_j^G(X)\lra \K_j^G(E)$ is an isomorphism.

\bth\label{equivkunneth} ({\it Equivariant Kunneth formula}) Let $X$
be a $G$-cellular scheme and let $Y$ be any $G$-scheme. Then the
canonical map $\varphi$ defined in (\ref{extalg}) is an isomorphism of
$K_{G}^0(Y)$-modules. Moreover, when $X$ and $Y$ are smooth then
$\varphi$ is an isomorphism of $K_{G}^0(Y)$-algebras. \eeth{\bf
  Proof:} Let $\alpha:X_{i+1}\hra X_i$ denote the embedding of the
closed subscheme and $\beta:U_i\lra X_i$ denote the immersion of the
open subscheme. For a $G$-scheme $Y$ consider the equivariant cellular
fibration
$$Y\times X_1\supseteq Y\times X_2\supseteq \cdots \supseteq Y\times
X_m=Y$$ over $Y$. Since $id_{Y}\times\alpha$ and
$id_{Y}\times \beta$ are maps of $G$-schemes, these induce
morphisms in higher $G$-equivariant $\K$-theory (\cite{Q}, \cite{Th})
$(id_{Y}\times\alpha)_*:\K_l^{G}(Y\times X_{i+1})\lra
\K^{G}_l(Y\times X_i)$ and
$(id_{Y}\times\beta)^*:\K_l^{G}(Y\times X_{i})\lra
\K_l^{G}(Y\times U_i)$ for each $l\geq 0$. We
have the following long exact sequence in higher $G$-equivariant
$\K$-theory \be\label{leseq1} \K_1^G(Y\times 
X_i)\stackrel{(id_{Y}\times \beta)^*}{\lra}
\K_1^{G}(Y\times U_i)\stackrel{\partial}{\lra} \K^G_0(Y\times
X_{i+1})\stackrel{(id_{Y}\times\alpha)_*}{\lra}
\K^{G}_0(Y\times X_i)\stackrel{(id_{Y}\times\beta)^*}{\lra}
\K^{G}_0(Y\times U_i)\lra 0\ee (see \cite[Lemma 5.5.1]{CG}).
Moreover, we have the pull back maps from
$\K
^{G}_l(Y)\lra \K^{G}_l(Y\times X_i)$ and the following commutative
triangle for every $l\geq 0$.
\[\begin{array}{lllll}
\K^{G}_l(Y\times X_i)&   &\stackrel{(id_{Y}\times\beta)^*}{\lra}&  &\K^{G}_l(Y\times U_i)\\
            &\nwarrow&    &\cong\nearrow &\\
             &        & \K^{G}_l(Y)&  &    
  \end{array}\] Note that the isomorphism in the above diagram is the
Thom isomorphism for the $G$-equivariant affine bundle $Y\times
U_i\lra Y$. It follows that the maps
$(id_{Y}\times\beta)^*:\K^{G}_l(Y\times X_i)\lra
\K^{G}_l(Y\times U_i)$ are
surjective. This in turn implies from (\ref{leseq1}) that the
connecting 
homomorphism $\partial$ is trivial and $(id_{Y}\times\alpha)_*$ is injective. Thus we have the following short exact sequence of $\K^0_{G}(Y)$-modules
\be\label{cfses1}0\ra \K_l^{G}(Y\times X_{i+1})\stackrel{{(id_Y\times \alpha)}_*}{\lra} \K_l^{G}(Y\times X_i)\stackrel{{(id_Y\times\beta)}^*}{\lra}
\K_l^{G}(Y\times U_i)\lra 0.\ee

When $Y=pt$ (\ref{cfses1}) reduces to the following short exact
sequence of $\R(G)$-modules 
\be\label{cfses2}0\ra \K_l^{G}(X_{i+1})\stackrel{{\alpha}_*}{\lra} \K_l^{G}(X_i)\stackrel{{\beta}^*}{\lra}
\K_l^{G}(U_i)\lra 0.\ee

We now claim that the map (\ref{map}) is an isomorphism for each
$1\leq i\leq m$. When $i=1$, this will imply that (\ref{extalg}) is an
isomorphism. We prove this by downward induction on $i$. This is
trivially true for $i=m$, since in this case $X_m=\emptyset$. Consider
the commutative diagram of $\R(G)$-modules
\be\label{diagramcf} \begin{array}{llllllllll}& &
  \K^{G}_0(Y)\bigotimes_{\R(G)} \K_0^{G}(X_{i+1}) &
  \stackrel{id_{Y}^*\bigotimes \alpha_*}{\ra} &
  \K^{G}_0(Y)\bigotimes_{\R(G)}
  \K_0^{G}(X_i)&\stackrel{id^*_{Y}\bigotimes\beta^*}{\ra}&
  \K_0^{G}(Y)\bigotimes_{\R(G)} \K_0^{G}(U_i)&\ra 0\\ &
  &\hspace{1cm}\da\phi_{i+1} & &\hspace{1cm}\da\phi_{i} &
  &\hspace{1cm}\da \psi_{i} & & \\ & 0\ra &~~~~\K_0^{G}(Y\times
  X_{i+1})&\stackrel{(id_{Y}\times
    \alpha)_*}{\ra}&~~~~\K_0^{G}(Y\times
  X_i)&\stackrel{(id_{Y}\times\beta)^*}{\ra}& ~~~~\K_0^{G}(Y\times
  U_i) \end{array}\ee where the bottom row is a part of (\ref{cfses1})
when $l=0$ and the top horizontal row is obtained from taking tensor
product of the exact sequence (\ref{cfses2}) when $l=0$ 
with the $\R(G)$-module $\K_0^{G}(Y)$ on the left.

Furthermore, the map $\psi_i$ is an isomorphism of $\R(G)$-modules for
every $1\leq i\leq m$. This can be seen from the following commuting
diagram:
\be\label{diagramthom} \begin{array}{llllllllll}
  \K^{G}_0(Y)\bigotimes_{\R(G)} \K_0^{G}(U_{i}) &
  \stackrel{\cong}{\lra} &
  \K^{G}_0(Y)\bigotimes_{\R(G)}
                                                \K_0^{G}(pt)\\
                         \hspace{1.2cm}\da\psi_{i} & & \hspace{1.2cm} \da\cong\\

                         \hspace{0.7cm} \K_0^{G}(Y\times
                         U_i) &\stackrel{\cong}{\lra}   &
                                                           \hspace{1cm}\K_0^G(Y)
                                                           \end{array}\ee

                                                           where the
                                                           top and the
                                                           bottom
                                                           isomorphisms
                                                           follow from
                                                           the Thom
                                                           isomorphisms
                                                           for the
                                                           $G$-equivariant
                                                           affine
                                                           bundles
                                                           $U_i\lra
                                                           pt$ and
                                                           $Y\times
                                                           U_i\lra Y$
                                                           respectively
                                                           and the
                                                           second
                                                           vertical
                                                           arrow is
                                                           the
                                                           canonical
                                                           isomorphism.

Therefore if we assume that $\phi_{i+1}$ is surjective, it follows by
diagram chase that $\phi_i$ is surjective. Note that the bottom
horizontal row of (\ref{diagramcf}) is left exact.
Hence if we assume that $\phi_{i+1}$ is injective, it follows again by
diagram chase that $\phi_i$ is injective. $\Box$

\bdefe\label{gcellrelative} Let $S$ be a $G$-scheme. By a relative {\em
  $G$-cellular scheme} we mean a $G$-scheme $X$ over $S$ equipped with
a $G$-stable relative algebraic cell decomposition. In other words
there is a filtration
\[X=X_1\supseteq X_2\supseteq\cdots\supseteq X_m\supseteq X_{m+1}=\emptyset\] where
each $X_i$ is a closed $G$-stable subscheme of $X$ over $S$ and
$X_i\setminus X_{i+1}=U_i$ is isomorphic to a $G$-equivariant affine
bundle which is a torsor under the vector group scheme associated to a
$G$-equivariant vector bundle over $S$ of rank $k_i$, for
$1\leq i\leq m$. Equivalently $X\lra S$ is a $G$-equivariant cellular
fibration in the sense of \cite{CG}.  \edefe

For $X$ and $Y$ schemes over $S$, we have a map of
$\K_{G}^0(S)$-modules
$$\varphi^S : \K^{G}_0(Y)\otimes_{\K_{G}^0(S)} \K^{G}_0(X)\lra
\K^{G}_0(Y\times_{S} X)$$ induced again by external tensor product of
$G$-equivariant coherent sheaves on $Y$ and $X$.

We state below the relative version of Theorem \ref{equivkunneth}.

\bth\label{equivkunneth1} Let $X$ be a relative $G$-cellular scheme
over $S$ and let $Y$ be any $G$-scheme over $S$. The map $\phi^S$ is
an isomorphism of $\K_{G}^0(Y)$-modules. Moreover, if $X$ and $Y$ are
smooth then $\phi^S$ is an isomorphism of
$\K_{G}^0(Y)$-algebras. \eeth {\bf Proof} The proof follows verbatim
as that of Theorem \ref{equivkunneth} by replacing everywhere
$Y\times X$ by the fibre product $Y\times_{S} X$ and $\R(G)$ by the
$\R(G)$-algebra $\K_{G}^0(S)$.  $\Box$

\brem\label {A-P} Since a $G$-cellular scheme is in particular
$T$-linear the proof of Theorem \ref{equivkunneth} also follows
directly from that of \cite[Proposition 6.4]{AP} by using the fact
that the $T$-equivariant Kunneth map is a faithfully flat extension of
the $G$-equivariant map, as in \cite[proof of Theorem A.3]{gon}. (Also
see Corollary \ref{corollary2} below where this argument is being
used.)  However the filtration argument we use in the proof of Theorem
\ref{equivkunneth} is better adapted to the setting of $G$-cellular
and $T$-cellular schemes, and allows an easy proof of Theorem
\ref{equivkunneth1}. Indeed the proof of Theorem \ref{equivkunneth1}
follows verbatim from that of \cite[Proposition 6.4]{AP} by replacing
everywhere $Y\times X$ by $Y\times_{S} X$. But it is good to see it
explicitly as above. \erem

\subsubsection{Applications of the Kunneth formula}

Let $X$ be a $T$-cellular scheme. In other words we have a
stratification \be\label{strat} X=X_1\supseteq
X_2\supseteq\cdots\supseteq X_m\supseteq X_{m+1}=\emptyset\ee by
$T$-stable closed subschemes such that
$F^{k_i}\cong U_i=X_i\setminus X_{i+1}$ is a $T$-representation.

Let $\cx:=G\times_{B} X$, where the $B$ action on $X$ is through the
canonical projection $B\lra T$. Then $\cx$ gets a $G$-scheme structure
via the natural action of $G$ on the left. We futher have a
$G$-equivariant stratification
$\cx=\cx_1\supseteq \cx_2\supseteq\cdots\supseteq \cx_m$ where
$\cx_i:=G\times_{B} X_i$ for $1\leq i\leq m$. Further,
$G\times_{B} U_i=\cx_i\setminus \cx_{i+1}$ is a $G$-equivariant vector
bundle over $G/B$. Thus $\cx$ has the structure of a relative
$G$-cellular scheme over $G/B$.

\bpropo\label{kunnethcf} Let $Y$ be any $G$-scheme. Then the canonical
map
\(\displaystyle\varphi: \K^{G}_0(Y)\bigotimes_{\R(G)} \K^{G}_0(\cx)
\lra \K^{G}_0(Y\times\cx)\) defined as in (\ref{extalg}) is an
isomorphism of $\K_{G}^0(Y)$-modules. Moreover, when $X$ and $Y$ are
smooth then $\varphi$ is an isomorphism of
$\K_{G}^0(Y)$-algebras. \epropo {\bf Proof:} Since $\cx$ is a relative
$G$-cellular scheme over $S=G/B$ by Theorem \ref{equivkunneth1} we
have the isomorphism \be\label{E1}
\K^{G}_0(Y_{S})\otimes_{\K^0_{G}(S)} \K^{G}_0(\cx)\cong
\K^{G}_0(Y_{S} \times_{S}\cx)\ee where $Y_{S}=Y \times S$ is the base
change of $Y$ to $S$. Since
$Y_{S}\times_{S} \cx= (Y\times S)\times_{S} \cx =Y\times \cx$, the
right hand side of (\ref{E1}) is isomorphic to
$\K^{G}_0(Y \times \cx)$.  By \cite[proof of Proposition 4.1,
p.30]{mer} we have the following isomorphism \be\label{E2}
\K^{G}_0(Y_{S})\cong \K^{G}_0(Y)\otimes_{\R(G)} \K^{G}_0(S)\ee since
$\K^G_0(S)=\R(B)=\R(T)$.  Furthermore, we have\be\label{E3}
\K^{G}_0(Y)\otimes_{\R(G)} \K^{G}_0(\cx)
\cong\K^{G}_0(Y)\otimes_{\R(G)}\K^{0}_G(S) \otimes_{\K^0_{G}(S)}
\K^{G}_0(\cx).\ee Note that the right hand side of (\ref{E3}) follows
by Theorem \ref{equivkunneth1} applied to $X=\cx$ and $Y=S$. Now, by
(\ref{E2}) and (\ref{E3}) we get that the left hand side of (\ref{E1})
is isomorphic to $\K^{G}_0(Y)\bigotimes_{\R(G)} \K^{G}_0(\cx)$. The
proposition now follows from (\ref{E1}). (Note here that
$\K_{G}^0(S)=\K^G_0(S)$ since $S$ is smooth.)  \hfill$\Box$

Let $\cb$ be a smooth base scheme and $p:\ce\lra \cb$ a
$G$-torsor. Let $\ce(X):=\ce\times_{G} X$ denote the associated bundle
with fibre a $G$-cellular scheme $X$ and projection
$\pi:\ce(X)\lra \cb$. We recall that $\ce$ is a smooth $G$-scheme and
$\ce(X)$ is a scheme.  Further, $\K_0(\ce(X))$ becomes a
$\K^0(\cb)$-module via pull back of classes of vector bundles under
$\pi^*$. Furthermore, we note that $\K^0(\cb)$ is an $\R(G)$-algebra
via the map which takes the isomorphism class of any
$G$-representation $\cv$ to the class in $\K^0(\cb)$ of the associated
vector bundle $\ce\times_{G} \cv$.

\bcor\label{maincor} We have the following isomorphism of
$\K^0(\cb)$-modules:
$$\K^0(\cb)\bigotimes_{R(G)} \K_0^{G}(X) \cong \K_0(\ce(X))$$ where the left
hand side has a canonical $\K^0(\cb)$-module structure by extension of
scalars to the $\R(G)$-algebra $\K^0(\cb)$. The above isomorphism is
an isomorphism of $\K^0(\cb)$-algebras if $X$ is a {\em smooth} $G$-cellular
scheme.\ecor

{\bf Proof:} Note that $X$ satisfies the hypothesis of Theorem
\ref{equivkunneth}. Further, since $G$ acts freely on $\ce$ as well as
on $\ce\times X$ diagonally, we have the
isomorphisms $$\K^0(\cb)=\K^0(\ce/G)=\K^0_{G}(\ce)$$
and $$\K_0(\ce(X))=\K_0((\ce\times X)/G)=\K^{G}_0(\ce\times X)$$ (see
\cite[Section 2.2]{mer} and \cite[Section 5.2.15]{CG}). Moreover, the
$\R(G)$-module structure on $\K^0_{G}(\ce)$ is via the pull back under
the structure morphism $\ce\lra \mbox{Spec}(F)$. Thus the class of a
$G$-representation $\cv$ pulls back to the class of the trivial bundle
$\ce\times \cv$ with the diagonal action of $G$ in $\K^0_{G}(\ce)$.
This further maps to the class in $\K^0(\cb)$ of the vector bundle
$\ce\times_{G} \cv$ over $\ce/G=\cb$. Note that since $\cb$ and hence
$\ce$ are smooth schemes
$K^0_{G}(\ce)=K^0(\cb)=K_0(\cb)=K^G_0(\ce)$. Also if in addition $X$
is a smooth $G$-scheme then $\ce\times X$ and hence $\ce(X)$ are
smooth schemes. The proof now follows readily from Theorem
\ref{equivkunneth}.  $\Box$

In the following corollary we show that the assertion of Corollary
\ref{maincor} holds under a weaker assumption that the $G$-scheme $X$
is $T$-cellular and not necessarily $G$-cellular. This is always true
if for instance we assume that $X$ is {\em smooth}, projective and has
only finitely many $T$-fixed points (see \cite{Bi1, Bi2, Bi3} or
\cite[Section 3.1, 3.2]{Br2}).

\bcor\label{corollary2} Let $X$ be $G$-scheme with a $T$-cellular
structure. We have the following isomorphism of $\K^0(\cb)$-modules
$$\K^0(\cb)\bigotimes_{\R(G)} \K^{G}_0(X) \cong \K_0(\ce(X)).$$ The
above isomorphism is an isomorphism of $\K^0(\cb)$-algebras if $X$ is
a {\em smooth} $T$-cellular $G$-scheme.\ecor {\bf Proof:} Since $X$ is
$T$-cellular we can apply Theorem \ref{equivkunneth} for the action of
$T$, taking $Y=\ce$. It follows that we have an isomorphism of
$\R(T)$-modules \be\label {equiv1}
\K_0^{T}(\ce)\bigotimes_{\R(T)}\K^{T}_0(X) \cong \K^T_0(\ce\times
X).\ee By \cite[Proposition 2.10]{mer} the isomorphism (\ref{equiv1})
can be rewritten as \be\label {equiv2}\K^{G}_0(\ce\times
G/B)\bigotimes_{\R(T)}\K^{G}_0(X\times G/B) \cong \K^G_0(\ce\times
X\times G/B).\ee Now, for a $G$-scheme $Y$ we have the following
canonical isomorphism of $\R(G)$-modules
\[ \R(T)\bigotimes_{\R(G)} \K^{G}_0(Y) \cong \K^{G}_0(Y\times G/B)\]
(see \cite[proof of Proposition 4.1, p.30]{mer}).  It follows that
(\ref{equiv2}) can be rewritten as \be\label{equiv3}
[\R(T)\bigotimes_{\R(G)}\K_0^{G}(\ce)]\bigotimes_{\R(T)}
[\R(T)\otimes_{\R(G)} \K^{G}_0(X)] \cong \R(T)\bigotimes_{\R(G)}
\K^{G}_0(\ce\times X).\ee Further, the left hand side of
(\ref{equiv3}) is isomorphic to
\(\displaystyle \R(T)\bigotimes_{\R(G)}
[\K^{G}_0(\ce)\bigotimes_{\R(G)} \K^G_0(X)]\). It follows that the
canonical map \be\label{map1}
\K^{G}_0(\ce)\bigotimes_{\R(G)}\K^{G}_0(X)\lra \K^{G}_0(\ce\times
X)\ee becomes an isomorphism after tensoring with $\R(T)$ which is a
free $\R(G)$-module of rank $|W|$ (see \cite[Proposition 1.22]{mer})
and hence a faithfully flat extension. (Also see \cite[proof of Theorem A3]{gon}.)  Therefore (\ref{map1}) must be an isomorphism. The
corollary follows by observing that since $\cb$ is smooth, we have
$\K^G_0(\ce)=\K_0(\cb)=\K^0(\cb)$. $\Box$

\subsection{Relative Localization theorem over $\cb=\ce/T$}

Throughout this subsection we let $\ce\lra \cb=\ce/T$ denote a $T$-torsor
and not a $G$-torsor.

Let $X$ be a {\em smooth} projective variety on which the torus
$T$ acts with finitely many $T$-fixed points.

We show below that $X$ admits {\it plus} and {\it minus} Bialynicki
Birula cell decomposition which are both filtrable (see
\cite{Bi1,Bi2,Bi3} and \cite[Section 3.1, 3.2]{Br2}). (The author is
grateful to Prof. M. Brion for the following explanation.)

Since $X$ is a smooth projective variety with $T$-action having only
finitely many fixed points $\{x_1,\ldots, x_m\}$, by a theorem of
Sumihiro (see \cite{sum}) there exists a $T$-equivariant embedding of
$X$ in a projective space $\mathbb{P}(V)$ where $V$ is a finite
dimensional $T$-module. We can write
$\displaystyle V=\bigoplus_{\chi\in X^*(T)} V_{\chi}$ which is a
decomposition of $V$ into $T$-eigenspaces. There are finitely many
characters $\chi's$ (say $\chi_1,\ldots, \chi_m$) such that $V_{\chi}$
is non zero. Thus we may find a one-parameter subgroup $\lambda$ of
$T$ such that the pairings $n_j:=\langle \chi_j,\lambda\rangle$ are
pairwise distinct for $i=1,\ldots, m$. We further reorder the indices
such that $n_1<\cdots<n_m$.

 We get an action of the multiplicative group $\mathbb{G}_m$ on
 $\mathbb{P}(V)$ via $\lambda$. The fixed points of this action are
 the projective subspaces $\mathbb{P}(V_{\chi_j})$ which are same
 as the fixed points for the action of $T$. Thus $\lambda$ is a {\it
   generic} one-parameter subgroup of $T$.

 Moreover, for a fixed point $x_i$ of $X$ we can define the plus cell
 \[X_+(x_i,\lambda):= \{ x\in X ~\mid ~ \lim_{t\rightarrow
     0}\lambda(t) \cdot x =x_i\}\] and the minus cell \[X_{-}(x,\lambda):=\{ x\in X ~\mid ~ \lim_{t\rightarrow
     \infty}\lambda(t) \cdot x =x_i\}.\]

 Thus for the fixed point $x_i=[v]$ for
 $\displaystyle v\in V_{\chi_{i}}$ it can be seen by direct
 computation that $X_{+}(x_i,\lambda)$ is the intersection of $X$ with
 the image in $\mathbb{P}(V)$ of
 $\displaystyle v+\bigoplus_{j>i} V_{\chi_j}$. Similarly
 $X_{-}(x_i,\lambda)$ is the intersection of $X$ with the image in
 $\mathbb{P}(V)$ of $\displaystyle v+\bigoplus_{j<i} V_{\chi_j}$.
 
 Note that the closure of the plus cell $X_{+}(x_i,\lambda)$ is
 contained in the union of this cell and
 $\displaystyle \mathbb{P}(\bigoplus_{j>i} V_{\chi_j})$. If we define
 the height of a $T$-fixed point
 $\displaystyle x\in \mathbb{P}(V_{\chi_i})$ as $h(x)=n_i$ then it can
 be seen that the closure of $X_{+}(x_i,\lambda)$ is contained in the
 union of this cell and the cells of fixed points with strictly bigger
 height and hence strictly bigger indices (by choosing a total order
 on $X^T$ compatible with the height function). Similarly, the closure
 of $X_{-}(x_i,\lambda)$ is contained in the union of this cell and
 the cells of fixed points with strictly smaller height.

 Thus we see that \begin{equation}\label{EQ1}\overline{X_{+}(x_i,\lambda)}\subseteq
   \bigcup_{j\geq i} X_{+}(x_j,\lambda)\end{equation}
 and \begin{equation}\label{EQ2}\overline{X_{-}(x_i,\lambda)} \subseteq \bigcup_{j\leq i}
   X_{-}(x_j,\lambda).\end{equation}

 In other words the  plus and  minus Bialynicki Birula cell
 decompositions of $X$ are filtrable.
 
 Moreover, it can be seen that the image in $\mathbb{P}(V)$ of
 $\displaystyle v+\bigoplus_{j>i} V_{\chi_j}$ and
 $\displaystyle v+\bigoplus_{j<i} V_{\chi_j}$ intersect transversally
 in $\mathbb{P}(V)$ at the $T$-fixed point $[v]$. It follows that
 $X_{+}(x_i,\lambda)$ and $X_{-}(x_i,\lambda)$ intersect transversally
 in $X$ at $x_i$ for $1\leq i\leq m$. Alternately this also follows
 from the decomposition
 $T_{x_i}X=(T_{x_i}X)_{-}\oplus (T_{x_i}X)_0\oplus (T_{x_i}X)_{+}$ of
 the tangent space of $X$ at $x_i$ into weight subspaces with
 negative, zero and positive weights under the $\mathbb{G}_m$ action
 via $\lambda$, and from the fact that that $(T_{x_i}X)_0$ is zero
 dimensional, $(T_{x_i}X)_{-}=T_{x_i} X_{-}(x_i,\lambda)$ and
 $(T_{x_i}X)_{+}=T_{x_i} X_{+}(x_i,\lambda)$.

 Let $U_i:=X_+(\lambda,x_i)$ and let
 $U_i':=X_{-}(\lambda,x_i)$. Further, let $V_i:=\bar{U_i}$ and
 $V_i':=\bar{U_i'}$. Then (\ref{EQ1}) and (\ref{EQ2}) can be rewritten
 as
 \begin{equation}\label{EQ3}V_i\subseteq \bigcup_{j\geq i} U_j\end{equation} and
\begin{equation}\label{EQ4}V'_i\subseteq \bigcup_{j\leq i} U'_j.\end{equation}

Furthermore, $(U_{i})^{T}=\{x_i\}$, $U_{i}\cong \mathbb{A}_{F}^{k_i}$,
$(U'_{i})^{T}=\{x_i\}$ and $U_{i}'\cong \mathbb{A}_{F}^{n-k_i}$
for every $1\leq i\leq m$. 

Then it follows from (\ref{EQ3}) and (\ref{EQ4}) that $V_i$ and $V_i'$
are closed $T$-invariant subvarieties of $X$ of dimensions $k_i$ and
$n-k_i$ respectively. Moreover, \be\label{intersection1} V_i~~
\mbox{and} ~~V_i' ~~\mbox{intersect transversally along}~~ x_i.\ee

It further follows from (\ref{EQ3}) and (\ref{EQ4}) that if $x_k\in V_i$
then $k\geq i$ and if $x_k\in V_j'$ then $k\leq j$. Thus if the
intersection of $V_i$ with $V'_j$ is nonempty then it must contain a
$T$-fixed point being a complete $T$-variety. Thus it follows that if
\be\label{intersection2} V_i\cap V_j'\neq \emptyset~~\implies ~~ i\leq
j.\ee

When $X$ is a smooth projective toric variety then $V_i$ and $V_i'$
are described explicitly in \cite[pages 102-104]{f}. When $X$ is the
flag variety then $V_i$ and $V_i'$ are the Schubert and the opposite
Schubert varieties respectively (see for example \cite[Section
1.2]{br}).

Further let \(\displaystyle X_{i}=\bigcup^m_{j=i} U_{j}\) and
\(\displaystyle X'_{i}=\bigcup^i_{j=1} U'_{j}\) so that
\[X_1=X\supseteq X_{2}\supseteq \cdots\supseteq X_m=\{x_m\}\] and
\[X'_m=X\supseteq X'_{m-1}\supseteq \cdots\supseteq X'_1=\{x_1\}\] are
the filtrations of $X$ by closed subvarieties given by the plus and
minus Bialynicki-Birula cell decompositions respectively.

We note that $\ce\times X\lra \ce$ is a $T$-equivariant cellular
fibration 
\[\ce\times X=\ce\times X_1\supseteq \ce\times X_2\supseteq
  \cdots\supseteq \ce\times X_m=\ce\times \{x_m\}.\] Further,
\(\displaystyle\ce\times U_i=\ce\times X_i\setminus \ce\times
X_{i+1}\) is isomorphic to a trivial $T$-equivariant vector bundle of
rank $k_i$
over $\ce\times x_i$. There is a dual $T$ -equivariant cellular
fibration structure on $\ce\times X\lra \ce$ \[\ce\times X=\ce\times X'_m\supseteq \ce\times X'_{m-1}\supseteq
  \cdots\supseteq \ce\times X'_{1}=\ce\times \{x_1\}\] where \(\displaystyle\ce\times U'_i=\ce\times X'_i\setminus \ce\times
X'_{i-1}\) is isomorphic to a trivial $T$-equivariant vector bundle
over $\ce\times x_i$ of rank $(n-k_i)$.

For $1\leq i\leq m$, consider
$\ce(V_i):=\ce\times_{T} V_i\subseteq \ce(X)$ and
$\ce(V'_i):=\ce\times_{T} V_i'\subseteq \ce(X)$ which are closed
subvarieties of $\ce(X)$ of codimensions $(n-k_i)$ and $k_i$
respectively. The following can be derived respectively from
(\ref{intersection1}) and (\ref{intersection2})

\be\label{intersection3} \ce(V_i)~~
\mbox{and} ~~\ce(V_i') ~~\mbox{intersect transversally along}~~ \ce(x_i).\ee

\be\label{intersection4} \ce(V_i)\cap
\ce(V_j')=\ce(V_i\cap V_j')\neq \emptyset~~\implies ~~ i\leq j.\ee

\bnotn Since $X$, $\cb$ and hence $\ce(X)$ are smooth, throughout this
section and the remaining part of the paper, henceforth we shall let
$\K_T(X):=\K^T_0(X)=\K^0_T(X)$, $\K(X):=\K_0(X)=\K^0(X)$,
$\K_T(\ce):=\K^T_0(\ce)=\K^0_T(\ce)$,
$\K_T(\ce\times X):=\K^T_0(\ce\times X)=\K^0_T(\ce \times X)$,
$\K(\cb):=\K_0(\cb)=\K^0(\cb)$ and
$\K(\ce(X)):=\K_0(\ce(X))=\K^0(\ce(X))$.\enotn

We shall fix a base point $b_0\in B$ and identify $X$ with the fibre
$\pi^{-1}(b_0)\subseteq \ce(X)$. We then have the restriction
homomorphism $\K(\ce(X))\lra K(X)$. We note here that the classes
$[\co_{\ce(V_i)}]\in \K(\ce(X))$ restrict to the classes $[\co_{V_i}]$
for $1\leq i\leq m$ which form a basis of $\K(X)$ as a free
$\mathbb{Z}$-module. The following is a version of Leray-Hirsch
theorem in this setting. (See \cite[p. 153, proof of Theorem 1.2
(iv)]{su} for the case when $X$ is a smooth projective toric variety.)

\bpropo\label{leray-hirsch} The classes
$[\co_{\ce(V_i)}]\in \K(\ce(X))$ for $1\leq i\leq m$ form a basis of
$\K(\ce(X))$ as a free $\K(\cb)$-module.  \epropo {\bf Proof:} The
classes of structure sheaves of the cell closures $[\co_{V_i}]_T$ form
a basis of $\K_T(X)$ as an $\R(T)$-module by \cite[Lemma 5.5.1(b)]{CG}. Now, under the isomorphism
$\K(\ce(X))\cong \K_T(\ce\times X) \cong \K(\cb)\otimes_{\R(T)}
\K_T(X)$ given by Corollary \ref{maincor} the classes
$[\co_{\ce(V_i)}]$ correspond to
$[\co_{\ce\times V_i}]_T\in \K_T(\ce\times X)$ and hence to the class
$1\otimes [\co_{V_i}]_T$ in $\K(\cb)\otimes_{\R(T)} \K_T(X)$ for
$1\leq i\leq m$. Since $1\otimes [\co_{V_i}]_T$ for $1\leq i\leq m$ is
a $\K(\cb)$-module basis of $\K(\cb)\otimes_{\R(T)} \K_T(X)$ the
proposition follows.$\Box$

In the rest of this section we shall assume that $X$ is a smooth
projective variety with an action of the torus $T$ having finitely many
$T$-fixed points as well as {\em finitely many $T$-invariant curves}.

Let \(\displaystyle \ce\times X^T\stackrel{\iota}{\hra} \ce\times X\)
denote the inclusion of $T$-subschemes. In this section we prove a
precise form of localization theorem for the $\K$-ring of the space
$\ce\times X$ which generalizes \cite[Theorem 1.3]{u2} to the relative
setting.

Let $C_{ij}$ denote the $T$-invariant irreducible curve in $X$ joining
the $T$-fixed points $x_i$ and $x_j$. Choose an isomorphism of
$C_{ij}$ with $\bp^1$ that sends $x_i$ to $0$ and $x_{j}$ to $\infty$
without loss of generality. Further, let $T$ act on
$C_{ij}\setminus x_j\cong \bp^1\setminus \infty$ via a character
$\chi_{ij}$. Let \(\displaystyle\mathcal{C}\) denote the finite collection
of invariant curves in $X$.

Let $\cy$ denote the subring of \(\displaystyle\prod_{k=1}^m \R(T)\)
consisting of $(y_k)$ such that $(1-e^{-\chi_{ij}})$ divides $y_i-y_j$ for
each $C_{ij}\in \mathcal{C}$. Clearly $\cy$ is an $\R(T)$-subalgebra of
\(\displaystyle\prod_{k=1}^m \R(T)\) where
\(\displaystyle \R(T)\hra \prod_{k=1}^m \R(T)\) is embedded diagonally.

Let $\cy_{ij}$ denote the subring of
\(\displaystyle\prod_{k=1}^m \R(T)\) consisting of $(y_k)$ satisfying
the condition that $(1-e^{-\chi_{ij}})$ divides $y_i-y_j$ corresponding to
$C_{ij}$ and $y_{k}\in \R(T)$ is arbitrary for $k\neq i,j$. Again
$\cy_{ij}$ is an $\R(T)$-subalgebra of
\(\displaystyle\prod_{k=1}^m \R(T)\) under the diagonal
embedding. Further, by definition
\be\label{intersection}\cy=\bigcap_{C_{ij}\in \mathcal{C}}
\cy_{ij}.\ee

Recall that $\K(\cb)=\K_{T}(\ce)$ has a canonical $\R(T)$-algebra
structure via the map that sends the class $[V]$ of a
$T$-representation to class $[\ce\times_{T} V]$ of the associated
vector bundle on $\cb$. In particular, $e^{\chi}\in \R(T)$ maps to the
class $[L_{\chi}]$ of the associated line bundle
$L_{\chi}:=\ce\times_{\R(T)} \mathbb{C}_{\chi}$.

For $1\leq i\leq m$, we have canonical sections
$s_i:\cb\lra \ce \times_{T} x_i\subseteq \ce\times_{T} X$ defined by
$s_i(b)=[e,x_i]$, where $e\in p^{-1}(b)$. Moreover,
$s_i$ and $s_j$ can be identified with the sections at $0$ and
$\infty$ of the $\bp^1$-bundle $\ce\times_{T} C_{ij}$ on $\cb$.

We have \be\label{resfp}\K_{T}(\ce\times
X^{T})=\K_{T}(\bigsqcup_{j=1}^m\ce \times x_j)\cong\prod_{j=1}^m
\K_{T}(\ce\times x_j) .\ee

Since $s_k$ maps $\cb$ isomorphically
onto $\ce(x_k):=\ce\times_{T} x_k$ with inverse
$\pi_k=\pi\mid_{\ce(x_k)}$ we have
$\K_{T}(\ce\times x_k)=\K(\ce(x_k))\cong
\K(\cb)$ for $1\leq k\leq m$. Thus
$s_k^*:\K(\ce(X))\lra \K(\cb)$ can be identified with the
composition of $\iota^*$ with the projection onto the direct factor
$\K_{T}(\ce\times x_k)\subseteq \K_{T}(\ce\times X^{T})$.

Further, $\iota^*$ can be identified with
\[(s_k^*): \K_{T}(\ce\times X)\lra \prod_{k=1}^m \K(\cb).\]

The ring $\K(\ce(X))$ gets the structure of $\K(\cb)$-algebra via pull
back by $\pi^*$, $\K(\cb)$-algebra structure on
$\displaystyle\prod_{k=1}^m \K(\cb)$ is via the diagonal inclusion and
$\iota^*$ is a morphism of $\K(\cb)$-algebras.

Let \(\displaystyle\cy':=\K_{T}(\ce)\bigotimes_{\R(T)} \cy\) and
\(\displaystyle\cy_{ij}':=\K_{T}(\ce)\bigotimes_{\R(T)} \cy_{ij}\) denote
respectively the extension of scalars of the $\R(T)$-algebras $\cy$ and
$\cy_{ij}$ to $\K_{T}(\ce)\cong \K(\cb)$. We can identify $\cy'$ with
\[ \{(y'_k)\in \prod_{k=1}^m \K(\cb)~\mid~1-[L_{\chi_{ij}}^{\vee}]~
  {\mbox {divides}}~
    y'_i-y'_j ~\forall~
C_{ij}\in \mathcal{C}\}.\]
Also $\cy'_{ij}$ can be identified with
\[ \{(y'_k)\in \prod_{k=1}^m \K(\cb)~\mid~1-[L_{\chi_{ij}}^{\vee}]~
  {\mbox {divides}}~
    y'_i-y'_j ~\mbox{corresponding~ to}
~C_{ij}\in \mathcal{C}~ \mbox{and}~y'_{k} ~\mbox {is arbitrary for} ~k\neq i,j\}.\] In particular,
$\cy'$ and $\cy'_{ij}$ are $\K(\cb)$-subalgebras of
\(\displaystyle\prod_{k=1}^m \K(\cb)\) under the diagonal
embedding. Furthermore, we note that
\be\label{intersection'}\cy'=\bigcap_{C_{ij}\in \mathcal{C}} \cy'_{ij}.\ee

\bth\label{relloc} Let $\ce\lra \cb$ be a $T$-torsor. Then the
restriction map 
\[\iota^*: \K_{T}(\ce\times X)\lra \K_{T}({\ce\times X}^{T})\cong
  \prod_{i=1}^m K(\cb)\] is injective and the image is isomorphic to
the $\K(\cb)$ subalgebra $\cy'$.  \eeth

{\bf Proof:} We first prove the injectivity of $\iota^*$.  By
Proposition \ref{leray-hirsch}, any element of $\K(\ce(X))$ is
uniquely expressible as a free $\K(\cb)$ module as follows:
\[ \sum_{i=1}^m \pi^*(b_i)\cdot [\co_{\ce(V_i)}]\] where
$b_i\in \K(\cb)$ for $1\leq i\leq m$. Let
\[\iota^*( \sum_{i=1}^m \pi^*(b_i)\cdot
  [\co_{\ce(V_i)}])=0.\] Let $k\geq 1$ be the least so that
$b_k\neq 0$.  Since $\iota^*$ is a $\K(\cb)$-algebra homomorphism this
implies that \be\label{eq} \iota^*([\co_{\ce(V_k')}])\cdot \sum_{i=k}^m
\pi^*(b_i)\cdot [\co_{\ce(V_i)}])=0.\ee

By (\ref{intersection4}), in $K(\ce(X))$ we have
\be\label{intersection5} [\co_{\ce(V_k')}]\cdot [\co_{\ce(V_i)}]=0\ee
whenever $i>k$. Thus from (\ref{intersection3}) and
(\ref{intersection5}), we get
$\displaystyle [\co_{\ce(V_k')}])\cdot \sum_{i=k}^m \pi^*(b_i)\cdot
[\co_{\ce(V_i)}]= \pi^*(b_k) \cdot [\co_{\ce(x_{k})}]$.
Thus (\ref{eq}) implies that 
\be\label{eq''} \iota^*(\pi^*(b_k) \cdot [\co_{\ce(x_{k})}])=0.\ee

Furthermore, $\iota^*$ is a $\K(\cb)$-algebra homomorphism where the
$\K(\cb)$-algebra structure on $\displaystyle \prod_{i=1}^m \K(\cb)$ is
via the diagonal embedding. Thus (\ref{eq''}) implies that
\be\label{eq'''} b_{k}\cdot \iota^* ([\co_{\ce(x_{k})}])=0.\ee

Since $\iota^*=(s_k^*)$,
(\ref{eq'''}) in particular implies that \be\label{eq'} b_{k}\cdot s_{k}^*([\co_{\ce(x_{k})}])=0.\ee

Now,
\[s_{k}^*([\co_{\ce(x_{k})}])= [\co_{B}]\] since $s_k$ is an isomorphism
of varieties from $\cb$ onto $\ce(x_k)$ whose inverse is
$\pi_k=\pi\mid_{\ce(x_k)}$. Thus from (\ref{eq'}) we get
$b_k\cdot [\co_{B}] =0$ in $\K(\cb)$. This implies that $b_k=0$ which
is a contradiction to our assumption. Thus we conclude that $\iota^*$
is injective.

Now, by Theorem \ref{equivkunneth} we have
$\K_{T}(\ce)\bigotimes_{\R(T)}\K_{T}(X)\cong \K_{T}(\ce\times X)$.
Moreover,
$\iota=id_{\ce}\times \iota': (\ce\times X)^{T}=\ce\times X^{T}\hra
\ce\times X$ where $\iota'$ denotes the inclusion $X^T\hra X$. Thus
$\iota^*$ can be identified with $id_{\K(\cb)}\otimes \iota'^*$.

More explicitly, by (\ref{resfp}) it follows that
\[\K_{T}(\ce\times X^{T})\cong \prod_{j=1}^m \K_{T}(\ce\times
  x_j)\cong \prod_{j=1}^m \K_{T}(\ce)\otimes_{\R(T)} \K_T(x_j) \cong
  \K_{T}(\ce)\bigotimes_{\R(T)}(\prod_{j=1}^m\K_{T}(x_j)) \cong
  \K_{T}(\ce)\bigotimes_{\R(T)}\K_{T}(X^T).\]

By \cite[Theorem
1.3]{u1}, the image of the restriction map
\(\displaystyle \iota'^*:\K_{T}(X)\hra \K_{T}(X^T)\cong \prod_{i=1}^m R(T)\)
is identified with the $\R(T)$-subalgebra $\cy$.  It follows that the
image of
\(\displaystyle\iota^*:\K_{T}(\ce)\bigotimes_{\R(T)}\K_{T}(X)\hra
\K_{T}(\ce)\bigotimes_{\R(T)} \K_{T}(X^T)\) can be identified with
$\cy'$. Hence the theorem. $\Box$

We have the following geometric interpretation of Theorem \ref{relloc}.
\bcor\label{locinterpretgeom}
We have a canonical embedding of $\K(\cb)$-algebras
\[\K(\ce(X))\stackrel{\iota^*}{\hra} \prod_{i=1}^m \K(\cb).\]
Furthermore, the image of $\iota^*$ is the intersection of the images
of
\[ \K(\ce(C_{ij}))\hra \K(\ce\times_{T} x_i)\times \K(\ce\times_{T}
  x_j)\hra \K(\ce\times_{T} X^T)\cong \prod_{i=1}^m \K(\cb).\] \ecor
{\bf Proof:} Recall that we can identify
\(\displaystyle \K_{T}(\ce\times X)\) with
\(\displaystyle \K(\ce(X))=(\ce\times X)/T\) and
\(\displaystyle \K_{T}(\ce \times x_i)=\K(\ce\times_{T} x_i)\)
canonically with the ring $\K(\cb=\ce/T)$ for every $1\leq i\leq
m$. Furthermore, Theorem \ref{relloc} applied for the principal
$T$-bundle $\ce\lra \cb$ and the smooth projective $T$-variety
$C_{ij}=\bp^1$ implies that
\(\displaystyle \K_{T}(\ce\times C_{ij})=\K(\ce(C_{ij}))\) embeds in
\(\displaystyle \K(\ce\times_{T} x_i)\times \K(\ce\times_{T}x_j)\hra
\K_{T}(\ce)\bigotimes_{R(T)} \K_{T}(X^T)\) and its image is isomorphic
to the subring \(\displaystyle \cy'_{ij}\). The proof now follows
readily from (\ref{intersection}) and Theorem \ref{relloc}. $\Box$

\subsection{Some further notations}

In this subsection and the remaining part of the paper we shall assume
that $G$ is an arbitrary reductive algebraic group and drop the
additional assumption of Subsections 1.1 and 1.2 that $\pi_1(G)$ is
torsion free.

Let $W$ denote the Weyl group and $\Phi$ denote the root system of
$(G,T)$. We further have the subset $\Phi^{+}$ of positive roots fixing $B\supseteq T$, and its
subset $\Delta=\{\alpha_1,\ldots, \alpha_r\}$ of simple roots where $r$
is the semisimple rank of $G$. For $\alpha\in\Delta$ we denote by
$s_{\alpha}$ the corresponding simple reflection. For any subset
$I\subset \Delta$, let $W_{I}$ denote the subgroup of $W$ generated by
all $s_{\alpha}$ for $\alpha\in I$. At the extremes we have
$W_{\emptyset}=\{1\}$ and $W_{\Delta}=W$.

Let $\Lambda:=X^*(T)$. Then $\R(T)$ (the representation ring of the
torus $T$) is isomorphic to the group algebra $\bz[{\Lambda}]$. Let
$e^{\lambda}$ denote the element of $\bz[{\Lambda}]=\R(T)$
corresponding to a weight $\lambda\in {\Lambda}$.  Then
$(e^{\lambda})_{\lambda\in {\Lambda}}$ is a basis of the $\bz$ module
$\bz[{\Lambda}]$. Further, since $W$ acts on $X^*(T)$, on
$\bz[{\Lambda}]$ we have the following natural action of $W$ given by
: $w(e^{\lambda})=e^{w(\lambda)}$ for each $w\in W$ and $\lambda\in
{\Lambda}$. Recall that we can identify $\R(G)$ with $\R(T)^{W}$ via
restriction to $T$, where $\R(T)^W$ denotes the subring of $R(T)$
invariant under the action of $W$ (see \cite[Example 1.19]{mer}).
 
Recall from \cite[Corollary 3.7]{Iv} that there exists an exact
sequence: \be\label{fc} 1\ra \cz\ra \tG:=\tC \times
G^{ss}{\stackrel{\pi}{\lr}}G \ra 1\ee where $\cz$ is a finite central
subgroup, $\tC$ is a torus and $G^{ss}$ is semisimple and
simply-connected.  The condition that $G^{ss}$ is simply connected
implies that $\tG$ is {\it factorial} (see \cite{mer}).

Now, from (\ref{fc}) it follows that $\tB:=\pi^{-1}(B)$ and
$\tT:=\pi^{-1}(T)$ are respectively a Borel subgroup and a maximal
torus of $\tG$. Further, by restricting the map $\pi$ to $\tT$ we get
the following exact sequence: \be\label{fc1} 1\ra \cz\ra \tT\ra T \ra
1. \ee Let $\tW$ and $\tphi$ denote respectively the Weyl group and
the root system of $(\tG,\tT)$.  Then by (\ref{fc}), it also follows
that $\tW=W$ and $\tphi=\Phi$.  Further we have \be\label{e1} \R(\tG)=
\R(\tC)\bigotimes \R(G^{ss})\ee and \be\label{e2} \R(\tT)\cong
\R(\tC)\bigotimes \R(T^{ss})\ee where $T^{ss}$ is the maximal torus
$\tT\cap G^{ss}$ of $G^{ss}$.

Recall we can identify $\R(\tG)$ with $\R(\tT)^{W}$ via restriction to
$\tT$, and further $\R(\tT)$ is a free $\R(\tG)$ module of rank $|W|$
(see \cite[Theorem 2.2]{st}). Moreover, since $G^{ss}$ is semi-simple
and simply connected, $\R(G^{ss})\cong \bz[x_1,\ldots, x_r]$ is a
polynomial ring on the fundamental representations (\cite[Example
1.20]{mer}). Hence $\R(\tG)=\R(\tC)\bigotimes \R(G^{ss})$ is the tensor
product of a polynomial ring and a Laurent polynomial ring, and hence
a regular ring of dimension $r+{dim}(\tC)={rank}(G)$ where $r$ is the
rank of $G^{ss}$.

For $X$ any $G$-scheme, we shall consider the $\tT$ and
$\tG$-equivariant $K$-theory of $X$ where we take the natural actions
of $\tT$ and $\tG$ on $X$ through the canonical surjections to $T$ and
$G$ respectively.

We consider $\bz$ as an $\R(\tG)$-module by the augmentation map
$\epsilon:\R(\tG)\ra \bz$ which maps any $\tG$-representation $V$ to
${dim}(V)$. Moreover, we have the natural restriction
homomorphisms $\K_{\tG}(X)\ra$ $\K_{\tT}(X)$ and $\K_{\tG}(X)\ra K(X)$
where $\K(X)$ denotes the ordinary Grothendieck ring of algebraic
vector bundles on $X$.  We then have the following isomorphisms (see
\cite[Proposition 4.1 and Theorem 4.2]{mer}).

\be\label{eqk1} \R(\tT)\bigotimes_{\R(\tG)}\K_{\tG}(X)\cong \K_{\tT}(X),\ee

 \be\label{eqk2}  \K_{\tG}(X)\cong \K_{\tT}(X)^{W},\ee

\be\label{eqk3} \bz\bigotimes_{\R(\tG)} \K_{\tG}(X)\cong \K(X).\ee

Let $\R(\tT)^{W_I}$ denote the invariant subring of the ring $\R(\tT)$
under the action of the subgroup $W_I$ of $W$ for every $I\subset
\Delta$. Thus in particular we have, $\R(\tT)^{W}=\R(\tG)$ and
$\R(\tT)^{\{1\}}=\R(\tT)$. Further, for every $I\subset \Delta$,
$\R(\tT)^{W_I}$ is a free module over $\R(\tG)=\R(\tT)^{W}$ of rank
$|W/W_{I}|$ (see \cite[Theorem 2.2] {st}). Indeed, \cite[Theorem
2.2]{st} which we apply here holds for $\R(T^{ss})$.  However, since
$W$ acts trivially on the central torus $\tC$ and hence trivially on
$\R(\tC)$ we have \be\label{e3}\R(\tT)^{W_I}=\R(\tC)\bigotimes
\R(T^{ss})^{W_{I}}\ee for every $I\subseteq \Delta$, and hence we
obtain the analogous statement for $\R(\tT)$.

Let $W^{I}$ denote the set of minimal length coset representatives of
the parabolic subgroup $W_I$ for every $I\subset \Delta$. Then
$$W^{I}:=\{w\in W \mid l(wv)=l(w)+l(v) ~\forall~ v\in W_I\}=\{w\in
W\mid w(\Phi_{I}^{+})\subset \Phi^{+}\}$$ where $\Phi_I$ is the root
system associated to $W_I$, where $I$ is the set of simple roots.
Recall (see \cite[p.19]{hum}) that we also have: 
$$W^{I}=\{w\in W\mid l(ws)>l(w)~for~all~s\in I\}.$$

Note that $J\subseteq I$ implies that $W^{\Delta\setminus
J}\subseteq W^{\Delta\setminus I}$.  Let
\be\label{1} C^{I}:=W^{\Delta\setminus I}\setminus (\bigcup_{J\subsetneq
I}W^{\Delta \setminus J}).\ee

Let $\alpha_1,\ldots,\alpha_r$ be an ordering of the set $\Delta$ of
simple roots and $\omega_1,\ldots,\omega_r$ denote respectively the
corresponding fundamental weights for the root system of
$(G^{ss},T^{ss})$.  Since $G^{ss}$ is simply connected, the
fundamental weights form a basis for $X^*(T^{ss})$ and hence for every
$\lambda\in X^*(T^{ss})$, $e^{\lambda}\in \R(T^{ss})$ is a Laurent
monomial in the elements $e^{\omega_i}:1\leq i\leq r$.

In \cite[Theorem 2.2]{st} Steinberg has defined a basis
$\{f_{v}^{^I}: v\in W^{I}\}$ of $\R(T^{ss})^{W_{I}}$ as an
$\R(T^{ss})^W$-module.  We recall here this definition: For $v\in
W^{I}$ let \be\label{st1} p_{v}:=\prod_{v^{-1}\alpha_{i}<0} e^{\omega_{i}}\in
\R(\tT).\ee Then \be\label{st2} f_v^{^{I}}:=\sum_{x\in W_{I}(v)\big{\backslash}
  W_{I}} x^{-1}v^{-1}p_{v}\ee where $W_{I}(v)$ denotes the stabilizer
of $v^{-1}p_v$ in $W_{I}$.

We shall also denote by $\{f_v^{^I}:v \in W^{I}\}$ the corresponding
basis of $\R(\tT)^{W_{I}}$ as an $\R(\tT)^{W}$-module where it is
understood that \be\label{notation} f_v^{^I}:=1\bigotimes f_v^{^I} \in
\R(\tC)\bigotimes \R(T^{ss})^{W_{I}}.\ee

\begin{notation}\label{aproposds} Whenever $v\in C^{I}$ we denote
$f_v^{^{\Delta\setminus I}}$ simply by $f_{v}$. We can drop the
superscript in the notation without any ambiguity since $\{C^{I}:
I\subseteq \Delta\}$ are disjoint. Therefore with the modified
notation \cite[Lemma 1.10]{u1} implies that: $\{f_{v}: ~v\in
W^{\Delta\setminus I}=\bigsqcup_{J\subseteq I} C^J\}$ form an
$\R(\tT)^{W}$-basis for $\R(\tT)^{W_{\Delta\setminus I}}$ for every
$I\subseteq \Delta$. Further, let \be\label{2}
\R(\tT)_{I}:=\bigoplus_{v\in C^I}\R(\tT)^{W}\cdot f_v .\ee \end{notation}

In $\R(T)$ let \be\label{multstr}f_{v}\cdot f_{v^{\prime}}=\sum_{J\subseteq
  (I\cup I^{\prime})} \sum_{w\in C^{J}}a^{w}_{v,v^{\prime}}\cdot f_{w}
\ee for certain elements $a^{w}_{v,v^{\prime}}\in \R(\tG)=\R(\tT)^{W}$
$~\forall~ v\in C^{I}$, $v^{\prime}\in C^{I^{\prime}}$ and $w\in
C^{J}$, $J\subseteq (I\cup I^{\prime})$.

\brem\label{simplyconnectedcover} In the next section, for a regular
compactification $X$ of $G$ we shall consider the action of
$\tG\times \tG$ on $X$ via its canonical surjection to $G\times G$ and
further consider $\K_{\tG\times \tG}(X)$ instead of
$\K_{G\times G}(X)$. This is in order to apply the results in
Subsections 1.1 and 1.2, since $\pi_1(\tG)$ is torsion free. Moreover,
this also enables us to use the Steinberg basis defined in Notation
\ref{aproposds} and its structure constants (\ref{multstr}) in the
description of the multiplicative structure of
$\K_{\tG\times \tG}(X)$.  \erem

\section{$\K$-theory of bundles with fibre regular compactifications of
  $G$}

In this section $X$ denotes a projective regular compactification of $G$.

Let $\bar T$ denote the closure of $T$ in $X$. It is known that for the
left action of $T$ (i.e. for the action of $T\times \{1\}$), $\bar{T}$
is a smooth projective toric variety. (see \cite{Br}).
Moreover, $X^{T\times T}$ is contained in the union
$X_c$ of all closed $G\times G$-orbits in $X$; moreover all such
orbits are isomorphic to $G/B^{-}\times G/B$.  

Let $\cf$ be the fan associated to $\bar T$ in $X_{*}(T)\otimes
\br$. Since $\bar{T}$ is complete, $\cf$ is a subdivision of
$X_{*}(T)\otimes \br$. There is a canonical action of
$\mbox{diag}(W):=\{(w_1,w_2)\in W\times W~ \mid ~w_1=w_2\}$ on
$\bar T$ induced from the conjugation action on $T$ which corresponds
to an action of $W$ on $\cf$. By \cite[Proposition A2]{Br}, it follows
that $\cf=W\cf_{+}$ where $\cf_{+}$ is the subdivision of the positive
Weyl chamber formed by the cones in $\cf$ contained in this
chamber. Therefore $\cf$ is a smooth subdivision of the fan associated
to the Weyl chambers, and $W$ acts on $\cf$ by reflection about the
Weyl chambers.  Let $\bar{T}^{+}$ denote the toric variety associated
to the fan ${\cf}_{+}$.  

Now, since $X$ is projective, the canonical morphism
$f : X \lra \bar{G_{ad}}$ is projective. Also, $\bar{T}^+$ is the
inverse image of $\mathbb{A}^r$ under $f$ and the restriction
$g : \bar{T}^+ \lra \mathbb{A}^r$ of $f$ is a projective morphism of
toric varieties. This implies that $\bar{T}^+$ is a semi-projective
$T$-toric variety.

Let ${\cf}(l)$ denote the set of maximal cones of $\cf$. Then we know
that ${\cf}_{+}(l)$ parameterizes the closed $G\times G$-orbits in
$X$.  Hence $X^{T\times T}$ is parametrized by
${\cf}_{+}(l)\times W\times W$(see \cite[Proposition A1 and A2]{Br}).

Recall by \cite[Theorem 2]{VV} and \cite[Theorem 2.1]{u1}, that
$\K_{\tT\times \tT}(X)$ embeds into $\K_{\tT\times \tT}(X_c)$, the latter being
a product of copies of the ring $\K_{\tT\times \tT}(G/B^{-}\times G/B)$.

Let $\cy$ denote $$(f_{\sigma,u,v})\in \prod_{\sigma\in \cf_+(l)} \prod_{u,v\in W\times W} \K_{\tT\times
  \tT}(x_{\sigma,u,v})=\K_{\tT\times \tT}(X^{\tT\times
  \tT})$$  satisfying the congruences:

\begin{enumerate}
\item[(i)] $f_{\sigma,us_{\alpha},vs_{\alpha}}\equiv f_{\sigma,
u,v}\pmod {(1-e^{-u(\alpha)}\bigotimes e^{-v(\alpha)})}$ whenever
$\alpha\in\Delta$ and the cone $\sigma\in{\cf}_{+}(l)$ has a facet
orthogonal to $\alpha$, and that

\item[(ii)] $f_{\sigma,u,v}\equiv f_{\sigma^{\prime}, u,v}\pmod
{(1-e^{-\chi})}$ whenever $\chi\in X^{*}(T)$ and the cones $\sigma$
and $\sigma^{\prime}\in{\cf}_{+}(l)$ have a common facet orthogonal to
$\chi$. 
 
\end{enumerate}
(In $(ii)$, $\chi$ is viewed as a character of $T\times T$ which is
trivial on $diag(T)$ and hence is a character of $T$.)

Then $\cy$ is an $\R(\tT)\otimes \R(\tT)$-subalgebra of
$\K_{\tT\times \tT}(X^{\tT\times \tT})$ (see \cite{u1}).

In this section we consider $\ce\lra \cb$ as a $\tG\times \tG$-torsor
over a scheme $\cb$. We shall consider the associated fibre bundle
\(\displaystyle \ce\times_{\tG\times \tG}X\) with fibre the regular
compactification $X$ of $G$ in view of Remark
\ref{simplyconnectedcover}.

The following proposition is the relative version of \cite[Theorem
2.1]{u1}.

\bth\label{locregularrel} Let $X$ be a projective regular
compactification of $G$ and let $\ce\lra \cb$ be a 
$\tG\times \tG$-torsor. The map \be\label{local}\prod_{\sigma\in
  \cf_{+}(l)}\iota_{\sigma}: \K_{\tT\times \tT}(\ce\times X)\lra
\prod_{\sigma\in \cf_{+}(l)} \K_{\tT\times \tT}(\ce\times G/B^-\times
G/B)\ee is injective and its image is
\(\displaystyle \K(\ce/\tT\times \tT)\bigotimes_{\R(\tT)\otimes \R(\tT)}
\cy\).  \eeth

{\bf Proof:} Consider the $\tT\times \tT$-torsor
$\ce\times \ce\lra \ce\times\ce/\tT\times \tT$. Consider the action of
$\tT\times \tT$ on $X$ by restriction. By Theorem \ref{relloc} we have
\be\label{step1}\iota:\K_{\tT\times \tT}(\ce\times X)\lra
\K_{\tT\times \tT}(\ce\times X^{T})\ee is injective. Further, by
applying Theorem \ref{equivkunneth} on either side of \ref{step1} we
get \be\label{step2} \iota:\K_{\tT\times \tT}(\ce\times X)\cong
\K(\ce/\tT\times \tT)\bigotimes_{\R(\tT)\bigotimes \R(\tT)} \K_{\tT\times \tT}(X) \lra
\K(\ce/\tT\times \tT)\bigotimes_{\R(\tT)\bigotimes \R(\tT)} \K_{\tT\times \tT}(X^T)\cong
\K_{\tT\times \tT}(\ce\times X^{T})\ee can be identified with the map
$id_{K(\ce/\tT\times \tT)}\bigotimes \iota'$, where $\iota'$ is the map
$\K_{\tT\times \tT}(X)\lra \K_{\tT\times \tT}(X^{T})$ induced by
restriction to the $T$-fixed points.  By \cite[Theorem 2.1]{u1} the
image of $\iota'$ lies in
\[\prod_{\sigma\in \cf_{+}(l)}\K_{\tT\times \tT}(G/B^-\times
G/B)\] and can be identified with $\cy$. Thus it follows that the image
of $\iota$ lies in
\[\K(\ce/\tT\times
\tT)\bigotimes_{\R(\tT)\bigotimes \R(\tT)} \prod_{\sigma\in
  \cf_{+}(l)}\K_{\tT\times \tT} (G/B^-\times G/B)\] and can be
identified with \(\displaystyle\K(\ce/\tT\times \tT)\bigotimes_{\R(\tT)\bigotimes \R(\tT)} \cy\). $\Box$

Let $\cz$ consist in all families $(f_{\sigma})_{\sigma\in
  \cf_{+}(l)}$ of elements of \(\displaystyle \R(\tT\times 1)\bigotimes \R(diag(\tT))\) such that
\begin{enumerate} 
\item[(i)] $(1,s_{\alpha})f_{\sigma}(u,v)\equiv f_{\sigma}(u,v) \pmod
{(1-e^{-\alpha(u)})}$ whenever $\alpha\in\Delta$ and
the cone $\sigma\in{\cf}_{+}(l)$ has a facet orthogonal to $\alpha$,
and that

\item[(ii)] $f_{\sigma}\equiv f_{\sigma^{\prime}}\pmod
{(1-e^{-\chi})}$ whenever $\chi\in X^{*}(T)$ and the cones $\sigma$ and
$\sigma^{\prime}\in{\cf}_{+}(l)$ have a common facet orthogonal to $\chi$.
 
\end{enumerate}
In particular, $\cz$ is $\R(\tG)\bigotimes \R(\tG)$-subalgebra of
\(\displaystyle\prod_{\sigma\in \cf_{+}(l)}
\R(\tT))\bigotimes \R(\tT)\).

The following is the relative version of \cite[Corollary 2.2]{u1}.

\bpropo\label{locregularinv} (i) We have a canonical inclusion
\be\label{inprodflagbun} \K(\ce(X))\hra
\prod_{\sigma\in \cf_{+}(l)}
\K(\ce/\tB^-\times \tB).\ee Here $\K(\ce/\tB^-\times \tB)$ is the
$\K$-ring of the bundle $\ce(G/B^-\times G/B)$ over $\cb$ with fibre
$G/B^-\times G/B$. 

(ii) The image of $\K(\ce(X))$ in the above inclusion is identified with $\K(\ce\times_{\tB\times \tB}\bar{T}^{+})$ which is the $\K$-ring of a toric bundle with fibre $\bar{T}^+$ over $\ce/B^-\times B=\ce(G/B^-\times G/B)$.

(iii) The ring $\K(\ce(X))$ is further isomorphic to
\(\displaystyle \K(\cb)\bigotimes_{\R(\tG)\otimes \R(\tG)} \cz.\) \epropo
{\bf Proof:} (i) By taking $W\times W$-invariants on either side of
(\ref{local}) in Proposition \ref{locregularrel} we get the inclusion
\be\label{locinv} [\K_{\tT\times \tT}(\ce\times X)]^{W\times W}\hra
\prod_{\sigma\in \cf_{+}(l)} [\K_{\tT\times \tT}(\ce\times G/B^-\times
G/B)]^{W\times W}.\ee Now, by applying \cite[Theorem 1.8]{u1} or
\cite[Proposition 4.1]{mer} on either side of (\ref{locinv}) we get:
\be\label{locinv1} \K_{\tG\times \tG}(\ce\times X)\hra \prod_{\sigma\in
  \cf_{+}(l)} \K_{\tG\times \tG}(\ce\times G/B^-\times G/B).\ee This is
further equivalent to \be\label{locinv2}\K(\ce\times_{\tG\times \tG}
X)\hra \prod_{\sigma\in \cf_{+}(l)} [\K(\ce\times_{\tG\times\tG}
\tG/\tB^-\times \tG/\tB)=\K(\ce/\tB^-\times \tB)]\ee and
(\ref{inprodflagbun}) follows.

(ii) Recall that we have a split exact sequence
\[ 1\lra \text{diag} ~\tT\lra \tT\times \tT\lra \tT\lra 1\] where the
second map is given by $(t_1,t_2)\mapsto t_1\cdot t_2^{-1}$ and the
splitting given by $t\mapsto (t,1)$. Thus we get canonical isomorphism
\be\label{splitchvar} \R(\text{diag}~ \tT)\bigotimes \R(\tT\times 1)\cong
\R(\tT\times \tT).\ee Using the change of variables coming from
(\ref{splitchvar}), \cite[Proposition 2.1]{u2} implies that the image
of $\K_{\tG\times \tG}(X)$ in
\(\displaystyle\prod_{\sigma\in \cf_+(l)} [\K_{\tG\times \tG}(G/B^-\times
G/B)=\R(\tT)\bigotimes \R(\tT)]\) can be identified with
$\K_{\tT}(\bar{T}^+)\bigotimes \R(\tT)$. Note that Corollary
\ref{corollary2} implies \be\label{kunnethreg} \K_{\tG\times
  \tG}(\ce\times X)\cong \K_{\tG\times \tG}(\ce)\bigotimes_{\R(\tG)\otimes
  \R(\tG)} \K_{\tG\times \tG}(X)\ee and \be\label{kunnethflag}
\K_{\tG\times \tG}(\ce\times G/B^-\times G/B)\cong \K_{\tG\times
  \tG}(\ce)\bigotimes_{\R(\tG)\bigotimes \R(\tG)} \K_{\tG\times
  \tG}(\tG/\tB^-\times \tG/\tB).\ee Thus under the inclusion
(\ref{locinv1}) the image of
\[\K_{\tG\times \tG}(\ce)\bigotimes_{\R(\tG)\otimes \R(\tG)}
\K_{\tG\times \tG}(X)\] in
\[\K_{\tG\times \tG}(\ce)\bigotimes_{\R(\tG)\otimes
  \R(\tG)}\prod_{\sigma\in \cf_{+}(l)}\R(\tT)\bigotimes \R(\tT)\] can be
identified with \be\label{image}\K_{\tG\times
  \tG}(\ce)\bigotimes_{\R(\tG)\bigotimes \R(\tG)} \K_{\tT}(\bar{T}^+)\bigotimes
\R(\tT).\ee By Theorem \ref{equivkunneth}, (\ref{image}) can further be
identified with
\[\K(\ce\times_{\tG\times \tG}(\tG\times\tG\times_{\tB^-\times
  \tB}\bar{T}^+\times pt))= \K(\ce\times_{\tB\times \tB} \bar{T}^+),\]
where $\tB^-\times \tB$ acts on $\bar{T}^+$ via the canonical
projection to $\tT\times 1$. 

(iii) Since \(\displaystyle\cz\cong \K_{\tG\times \tG}(X)\) by \cite[Proposition
2.5]{u1} and \(\displaystyle \K(\cb)\cong \K_{\tG\times \tG}(\ce)\), the claim readily
follows from (\ref{kunnethreg}).  $\Box$

\subsection{First description of $\K(\ce(X))$}

Recall that
$$\ce/(\tB^-\times \tB)=\ce\times_{\tG\times \tG}
(\tG\times\tG)/(\tB^-\times\tB) $$ is a bundle with fibre
$\tG/\tB^-\times\tG/\tB$ over $\cb$. Thus
\(\displaystyle \K(\ce(\tG/\tB^-\times\tG/\tB))\) gets a
$R(\tB^-)\times R(\tB)$-module structure by sending a representation
$\cv\otimes \cw$ of $\tB^-\times \tB$ to the associated vector bundle
\(\displaystyle \ce\times_{\tB^-\times \tB}(\cv\otimes \cw)\) on
\(\displaystyle\ce/(\tB^-\times \tB)\). Since \(\displaystyle
\tG\times\tG\times_{\tB^-\times \tB}\cv\otimes \cw\) is a $\tG\times
\tG$-linearized vector bundle on the space \(\displaystyle \tG/\tB^-\times\tG/\tB\). This is also the associated bundle \(\displaystyle\ce\times_{\tG\times\tG}(\tG\times\tG\times_{\tB^-\times \tB}\cv\otimes \cw)\). Let
\(\displaystyle \K_{\tB^-\times \tB}(\bar{T})\) denote the $\tB^-\times \tB$-equivariant $\K$-ring of $\bar{T}$ where we take the natural action of $\tB^-\times \tB$ on $\bar{T}$ via the canonical projection to $T\times T$.

We now prove the first main theorem of this section.

\bth\label{first} Let $\ce\lra \cb$ be a $\tG\times
\tG$-torsor. Consider the associated bundle
$\ce(X):=\ce\times_{\tG\times\tG} X$ with fibre the regular
compactification $X$ of $G$ over $\cb$. Here again the action on
$\tG\times \tG$ on $X$ is via the natural projection to $G\times
G$. The ring $\K(\ce(X))$ is isomorphic to the ring
$\K(\ce\times_{\tB^-\times \tB} \bar{T})^{\mbox{diag} (W)}=$
\be\label{expflagtorbun}[\K(\ce((\tG\times
\tG)/(\tB^-\times\tB)))\bigotimes_{\R(\tB^-)\bigotimes
  \R(\tB)}\K_{\tB^-\times \tB}(\bar{T})]^{\mbox{diag}(W)}\ee as a
$\K(\cb)$-module.  \eeth {\bf Proof:} By Corollary \ref{corollary2} we
have \be\label{eq1} \K(\ce(X))\cong \K(\cb)\bigotimes _{\R(\tG\times
  \tG)}\K_{\tG\times \tG}(X).\ee By \cite[Corollary 2.3]{u1}
(\ref{eq1}) implies \be\label{eq2} \K(\ce(X))\cong\K(\cb)\bigotimes
_{\R(\tG\times \tG)}\K_{\tT\times \tT}(\bar{T})^{\mbox{diag}(W)}.\ee
By \cite[Corollary 2.15]{mer} this can further be rewritten as
\be\label{eq3} \K(\ce(X))\cong \K(\cb)\bigotimes _{\R(\tG\times
  \tG)}\K_{\tB^-\times \tB}(\bar{T})^{\mbox{diag}(W)}\ee and
\be\label{eq4} \K(\ce(X))\cong \K(\cb)\bigotimes _{\R(\tG\times
  \tG)}\K_{\tG\times \tG}(\tG\times \tG\times_{\tB^-\times
  \tB}\bar{T})^{\mbox{diag}(W)}.\ee Since $\mbox{diag}(W)$ acts
trivially on $\cb$ and hence on $\K(\cb)$, by Proposition
\ref{kunnethcf} it follows that the right hand side of (\ref{eq4}) is
isomorphic to
$\K(\ce\times_{\tG\times\tG}(\tG\times\tG\times_{\tB^-\times \tB}
\bar{T}))^{\mbox{diag}(W)}$. This reduces to
$\K(\ce\times_{\tB^-\times\tB} \bar{T})^{\mbox{diag}(W)}$. Now,
(\ref{expflagtorbun}) follows by applying Corollary \ref{maincor} to
the $\tB^-\times \tB$-torsor $\ce\lra \ce/(\tB^-\times \tB)$ and the
associated $\bar{T}$-bundle.  $\Box$

\subsection{Second description of $\K(\ce(X))$}
We first set up some notations necessary to state the main theorem.

Consider the ring \be\label{ktorbunflag} \ck:=\K(\ce\times_{\tB\times
  \tG}\bar{T}^{+})\ee where  $\tB\times \tG$-acts on $\bar{T}^+$ via
the canonical projection $\tB\times \tG\lra \tT\times 1$. The ring (\ref{ktorbunflag}) is the $\K$-ring of a $\bar{T}^+$-bundle over the flag bundle $\ce/\tB\times \tG=\ce\times_{\tG\times \tG}(\tG/\tB \times pt)$ over $\cb$. Since $\bar{T}^+$ is a semi-projective  toric variety, by \cite[Theorem 4.1]{u2} the ring $\ck$ gets a $\K(\ce/\tB\times \tG)$-algebra structure.

\be\label{basis} \mbox{Let} ~\bar{f_v}:=1\bigotimes (1\bigotimes f_v)\in \K(\cb)\bigotimes_{\R(\tG)\bigotimes \R(\tG)} \R(\tG)\bigotimes \R(\tT)= \K(\ce/{\tG\times \tB})\ee where $f_v\in \R(\tT)=\K_{\tG}(\tG/\tB)$ is as in Notation \ref{aproposds}.  Note that  $\ce/{\tG\times \tB}=\ce\times_{\tG\times\tG} (pt\times\tG/\tB)$ is a flag bundle over $\cb=\ce/{\tG\times \tG}$. 

\be\label{eulerclass} \mbox{Let}~ \lambda_{I}:=1\bigotimes (\mu_{I}\bigotimes 1)\in \K(\cb)\bigotimes_ {\R(\tG)\bigotimes \R(\tG)} \R(\tT)\bigotimes \R(\tG)=\K(\ce/\tB\times \tG)\ee where \(\displaystyle\mu_{I}:=\prod_{\alpha\in I} (1-e^{-\alpha})\in \R(\tT)\) for $I\subset \Delta$. 

\be\label{coeffstein}  \mbox{Let} ~ c^w_{v,v'}:= 1\bigotimes (1\bigotimes a^w_{v,v'})\in \K(\cb)\bigotimes_ {\R(\tG)\bigotimes \R(\tG)} \R(\tG)\bigotimes \R(\tG)= \K(\cb)\ee where $a^w_{v,v'}\in \R(\tG)$ is as in (\ref{multstr}). 

Let
\be\label{ktorbunflag'} {\ck'}:=\K(\ce\times_{\tB \times \tB}
\bar{T}^{+}).\ee Here $\tB\times \tB$-acts on $\bar{T}$ via the
canonical projection to $\tT\times 1$. Then $\ck'$ is the $\K$-ring of a $\bar{T}^+$-bundle over the bundle $\ce/B^-\times B$ having fibre $G/B^-\times G/B$ over $\cb$. Again since $\bar{T}^+$ is a semi-projective  toric variety, by \cite[Theorem 4.1]{u2} the ring ${\ck'}$ gets a $\K(\ce/\tB\times \tB)$-algebra structure.

Further, we note that $\ce/\tB\times \tB$ is a flag bundle over $\ce/\tB\times \tG$ with fibre the flag variety $pt\times \tG/\tB$. Moreover, $\ce\times_{\tB\times\tB}\bar{T}^+$ is the pull back of $\ce\times_{\tB\times \tG} \bar{T}^+$ to $\ce/\tB\times \tB$. Thus the canonical inclusion \be \K(\ce/\tB\times\tG)\hra \K(\ce/\tB\times \tB)\ee is the restriction of $\ck\hra {\ck}'$. 

Moreover, $\bar{f_v}$, $\lambda_{I}$ and $c^w_{u,v}$ lie in $\K(\ce/\tB\times \tB)$ via pull back from $\K(\ce/\tG\times \tB)$, $\K(\ce/\tB\times \tG)$ and $\K(\cb)$ respectively.

We now prove the second main theorem of this section. This is the
relative version of \cite[Theorem 3.8]{u1}.

\bth\label{second}  \begin{enumerate}
\item We have the following isomorphism as submodules of ${\ck}'$: \be\label{isosub} \K(\ce(X))\cong \bigoplus_{v\in \Delta} \ck \cdot \bar{f_v}.\ee In particular, the ring $\K(\ce(X))$ gets a canonical structure of a $\ck$-module of rank $|W|$.  

\item Furthermore, (\ref{isosub}) is an isomorphism of $\ck$-algebras where any two basis elements $\bar{f_v}$ and $\bar{f_{v'}}$ multiply in $\ck'$ as follows \be\label{multbasis} \bar{f_v}\cdot \bar{f_{v'}}:=\sum_{J\subseteq I\cup I'}\sum_{w\in C^J} (\lambda_{I\cap I'}\cdot \lambda_{(I\cup I')\setminus J})\cdot c^w_{v,v'}\cdot \bar{f_w}.\ee

\end{enumerate}
\eeth
{\bf Proof:} \begin{enumerate}
\item Note that the isomorphism in \cite[Theorem 2.2 (i)]{u2} is as $\R(\tG)\bigotimes \R(\tG)$-algebras. Thus by base changing to
  the $\R(\tG)\bigotimes \R(\tG)$-module $\K(\cb)$ on either side we
  get the following isomorphism of $\K(\cb)$-algebras
  $$ \K(\cb)\bigotimes_{\R(\tG)\bigotimes \R(\tG)} \K_{\tG\times \tG}(X)
  \cong \bigoplus_{I\subseteq \Delta}\K(\cb)\bigotimes_ {\R(\tG)\otimes
    \R(\tG)} \K_{\tT}(\bar{T}^+)\bigotimes \R(\tT)_{I}.$$ By Corollary
  \ref{maincor} this can be rewritten as
  $$\K(\ce(X))\cong\bigoplus_{I\subseteq \Delta}\bigoplus_{v\in C^I}
  \K(\cb)\bigotimes_{\R(\tG)\bigotimes \R(\tG)} \K_{\tT}(\bar{T}^+)\bigotimes
  \R(\tG)\cdot f_v .$$ Now, $\R(\tG)\bigotimes \R(\tG)$ acts on
  $\K_{\tT}(\bar{T}^+)$ and $\R(\tG)\cdot f_v$ via the first and second
  projections respectively. Thus (\ref{ktorbunflag}) and (\ref{basis})
  together imply (\ref{isosub}). Note that $\ck\cdot \bar{f_v}$ is a
  $\ck$-submodule of $\ck'$ for every $v\in C^{I}$ and
  $I\subseteq \Delta$. Furthermore, the direct sum decomposition
  (\ref{isosub}) gives $\K(\ce(X))$ a structure of a free $\ck$-module
  of rank $|W|$. Also by Proposition \ref{locregularinv} (ii),
  (\ref{isosub}) is an equality of $\ck$-submodules of $\ck'$.

\item We observe that
  \[\bar{f_v}\cdot \bar{f_{v'}}=[1\bigotimes (1\bigotimes f_v)]\cdot
    [1\bigotimes (1\bigotimes f_{v'})]=1\bigotimes [(1\bigotimes
    f_v)\cdot (1\bigotimes f_{v'})].\] Now, \cite[Theorem 2.2
  (ii)]{u2} implies that $\bar{f_v}\cdot \bar{f_{v'}}=$
  \[1\bigotimes \sum_{J\subseteq (I\cup I^{\prime})}\sum_{w\in
      C^{J}}(\mu_{ I\cap I^{\prime}}\cdot \mu_ {I\cup
      I^{\prime}\setminus J}\bigotimes a^{w}_{v,v^{\prime}})\cdot
    (1\bigotimes f_{w}).\] This can further be written
  as
  \[1\bigotimes \sum_{J\subseteq (I\cup I^{\prime})}\sum_{w\in
      C^{J}}(\mu_{ I\cap I^{\prime}}\cdot \mu_ {I\cup
      I^{\prime}\setminus J}\bigotimes 1)\cdot (1\bigotimes
    a^{w}_{v,v^{\prime}})\cdot (1\bigotimes f_{w}).\] The equality
  (\ref{multbasis}) now follows by applying (\ref{basis}),
  (\ref{eulerclass}) and (\ref{coeffstein}) successively.
\end{enumerate} $\Box$

\beg
Let $X_0$ denote the wonderful compactification of
$PGL(n,\mathbb{C})$.  Now, regarding $\mathbb{P}^n(\mathbb{C})$ as
$PGL(n+1,\mathbb{C})/Q$ where $Q$ denotes the maximal parabolic in
$PGL(n+1,\mathbb{C})$ which fixes the one dimensional subspace of
$\mathbb{C}^{n+1}$ generated by the coordinate vector $e_1$ we can
consider the principal $Q\times Q$-bundle
\[PGL(n+1,\mathbb{C})\times PGL(n+1,\mathbb{C})\lra
  \mathbb{P}^n(\mathbb{C})\times \mathbb{P}^n(\mathbb{C})\cong
  PGL(n+1,\mathbb{C})\times PGL(n+1,\mathbb{C})/Q\times Q.\] The Levi
subgroup $L_{Q}$ of $Q$ is identified with $GL(n,\mathbb{C})$ and the
adjoint group $L_{Q}/C(L_{Q})$ is $PGL(n,\mathbb{C})$. We can
construct the associated bundle
\[X_1:=(PGL(n+1,\mathbb{C})\times PGL(n+1,\mathbb{C}))\times_{Q\times Q}
X_0\] where $Q\times Q$ acts on
$X_0$ via its projection to
$L_{Q}/C(L_{Q})\times L_{Q}/C(L_{Q})$. Since $Q/L_{Q}=Q_u$ is the
unipotent radical of $Q$, we have
\[\K(X_1)\cong \K(PGL(n+1,\mathbb{C})\times PGL(n+1,\mathbb{C})
  \times_{L_{Q}\times L_{Q}} X_0) \]  and
\[ \K(PGL(n+1,\mathbb{C})\times PGL(n+1,\mathbb{C})/Q\times Q)\cong
  \K( PGL(n+1,\mathbb{C})\times PGL(n+1,\mathbb{C}) /L_{Q}\times
  L_{Q})\cong \K(\mathbb{P}^n(\mathbb{C})\times
  \mathbb{P}^n(\mathbb{C}))\] (see \cite[5.2.18]{CG}).  Thus by
Theorem \ref{first} we have
$\K(X_1)\cong \K(\mathbb{P}^n(\mathbb{C})\times
\mathbb{P}^n(\mathbb{C}))\otimes_{\R(L_{Q})\otimes \R(L_{Q})}
\K_{L_{Q}\times L_{Q}}(X_0)$. Here
$\R(L_{Q}\times L_{Q})\cong \R(T\times T)^{S_n\times S_n}$ (see
Section 1.3 above) and
$\K_{L_{Q}\times L_{Q}}(X_0)\cong \K_{T\times T}(Y_0)^{S_n}$ where
$T$ the maximal torus of $GL(n,\mathbb{C})$ consisting of the diagonal
matrices acts canonically on the toric variety $Y_0$ associated to the
Weyl chambers of $A_{n-1}$-type. Thus we have
\[\K(X_1)=\K(\mathbb{P}^n(\mathbb{C})\times
  \mathbb{P}^n(\mathbb{C}))\otimes_{\R(T\times T)^{S_n\times S_n}}
    \K_{T\times T}(Y_0)^{S_n}.\]\eeg

\section{$K$-theory of toric bundles and flag bundles}
In this section we retrieve known results on $K$-theory of toric
bundles \cite[Theorem 1.2(iv)]{su} and flag bundles.

Let $X$ be a smooth $T$-cellular toric variety associated to a fan
$\Sigma$ in the lattice $N=\bz^n$. Let
$\Sigma(1)=\{\rho_1,\ldots,\rho_d \}$ denote the edges and
$v_1,\ldots,v_d$ primitive lattice points along the edges. Let
$M=Hom(N,\bz)$ be the dual lattice.

Let $\ce\lra \cb$ be a principal $T$-bundle. Let $\ce(X)$ denote the
associated toric bundle $\ce\times_{T} X$.

\bpropo\label{toric}
 Then $\K(\ce(X))$ has the following presentation as $\K(\cb)$-algebra:
\[\K(\cb)[ x_1,\ldots,x_d]/\ci \]where $\ci$ is the ideal generated by the following two types of relations:

\be\label{rel1} x_{i_1}\cdots x_{i_k}~~\mid~~\langle v_{i_1},\ldots, v_{i_k}\rangle\notin \Sigma.\ee

\be\label{rel2} \prod_{i\mid \langle u,v_i\rangle\geq 0}
(1-x_i)^{\langle u,v_i\rangle}- [L_u]\prod_{i\mid \langle
  u,v_i\rangle\leq 0} (1-x_i)^{-\langle u,v_i\rangle}~\forall~~u\in
M\ee \epropo {\bf Proof:} Since $X$ is $T$-cellular it satisfies the
hypothesis of the Theorem \ref{equivkunneth} and Corollary
\ref{maincor}. Hence by Corollary \ref{maincor},
\(\displaystyle \K(\ce(X))=\K(\cb)\bigotimes_{R(T)}\K_{T}(X)\) where
the extension of scalars to $\K(\cb)$ is obtained by sending
$e^u\in \R(T)$ to $[L_u]$ for every $u\in M=Hom(T,\bc^*)$. Now the
theorem follows readily from the presentation of the ring $\K_{T}(X)$
as an $\R(T)$-algebra described in \cite[Theorem 6.4]{VV} . $\Box$

Let $\ce\lra \cb$ be a principal $T$-bundle and let $\ce(G/B)$
denote the associated flag bundle $\ce\times_{T} G/B$.

\bpropo\label{flag2} The ring $\K(\ce(G/B))$ has the following
presentation \be\label{flagpres} \frac{\K(\cb)\bigotimes
  \R(T)}{\ci}\ee where $\ci$ is the ideal generated by the
relations
\[[\ce\times_{G} \cv]\otimes 1-1\otimes [\cv]~for~every~[\cv]\in
  \R(G)=\R(T)^{W}.\] \epropo {\bf Proof:} Now $G/B$ is a projective
variety with a $T$-cellular structure given by the Bruhat
decomposition. Therefore by Corollary \ref{corollary2} we have
\(\displaystyle \K(\ce(G/B))=\K(\cb)\bigotimes_{\R(G)}
\K_{G}(G/B)\). The theorem now follows by using the fact that
$\K_{G}(G/B)=\R(T)$ and the $\R(G)$ algebra structure on $\R(T)$ and
$\K(\cb)$ (see \cite[(2)]{pr}).$\Box$.

We further derive the following description of $\K(\ce(G/B))$.

\bpropo\label{localflagbun} The ring $\K(\ce(G/B))$ is isomorphic to
the subring of $\K(\cb)^{|W|}$ consisting of tuples $(f_w)_{w\in W}$
satisfying the condition that $f_w-f_{w\cdot s_{\alpha}}$ is divisible
by $1-[L^{\vee}_{\alpha}]$ for every $w\in W$ and $\alpha\in \Phi^+$.
\epropo {\bf Proof:} Recall that the $T$-fixed point set in $G/B$ can
be identified with $W$. Thus from the description of the $T$-invariant
curves in $G/B$ (see \cite{C} or \cite[Section 6.5]{Br2}) and from
\cite[Theorem 1.3]{u1} it follows that $\K_T(G/B)$ can be described
as the $R(T)$-subalgebra of $\R(T)^{|W|}$ consisting of tuples $(f_w)$
such that $f_w-f_{s_{\alpha}w}$ is divisible by $1-e^{-\alpha}$ for
$w\in W$ and $\alpha\in \Phi^+$. The proposition now follows
readily from Theorem \ref{relloc}. $\Box$

\end{document}